\documentclass[12pt]{amsart}

\usepackage{graphicx}
\usepackage{amsfonts}
\usepackage{epsf}
\usepackage{amssymb}
\usepackage{amsmath}
\usepackage{color}

\newtheorem{theorem}{Theorem}[section]
\newtheorem{proposition}[theorem]{Proposition}
\newtheorem{lemma}[theorem]{Lemma}
\newtheorem{definition}[theorem]{Definition}

\newtheorem{corollary}[theorem]{Corollary}
 
\newtheorem{remark}[theorem]{Remark} 
\newtheorem{example}[theorem]{Example}

\newcommand{\eps}{\varepsilon}
\newcommand{\qq}{\mathbb{Q}}
\newcommand{\dotf}{\dot{\mathbb{F}}}
\newcommand{\dotz}{\dot{\mathbb{Z}}}
\newcommand{\lr}{\langle \cdot , \cdot \rangle}

\begin{document}

\title{The rational Witt class and the unknotting number of a knot}
\thanks{The author was partially supported by NSF grant DMS 0709625.}
\begin{abstract}
We use the rational Witt class of a knot in $S^3$ as a tool for addressing questions about its unknotting number. We apply these tools to several low crossing knots (151 knots with 11 crossing and 100 knots with 12 crossings) and to the family of $n$-stranded pretzel knots for various values of $n\ge 3$. In many cases we obtain new lower bounds and in some cases explicit values for their unknotting numbers. Our results are mainly concerned with unknotting number one but we also address, somewhat more marginally, the case of higher unknotting numbers. 
\end{abstract}
\author{Stanislav Jabuka}
\email{jabuka@unr.edu}
\address{Department of Mathematics and Statistics, University of Nevada, Reno NV 89557, USA.}
\maketitle
\section{Introduction}
\subsection{Preliminaries and statement of results}
The {\em unknotting number $u(K)$} of a knot $K$ in the $3$-sphere is the minimum number of crossing changes, in any regular projection of $K$, that renders it unknotted. While $u(K)$ is easy to define, computing it in practice is often unwieldy. Some of the lower bounds for $u(K)$ come from the Tristram-Levine signatures $\sigma _\omega (K)$, $\omega \in S^1$ (see Definition \ref{definition-of-tristram-levine-signature}) and bound $u(K)$ as\footnote{We indicate a simple proof of this bound at the end of Section \ref{section-phi-of-k-under-a-crossing-change}. The usual knot signature $\sigma (K)$ agrees with $\sigma _{-1}(K)$. } 
\begin{equation} \label{lower-bounds-on-uK}
|\sigma _{\omega}(K)|\le 2 \, u(K) \quad \quad \quad \forall  \, \omega \in S^1
\end{equation}
On the other hand, upper bounds for $u(K)$ are most easily found from explicit unknottings of $K$. It is when the upper and lower bounds are disparate, that $u(K)$ is difficult to determine.  

The last three decades have furnished an impressive array of tools for studying unknotting numbers, tools   stemming from varied sources such as gauge theory \cite{cochran-lickorish, ozsvath-szabo, owens}, polynomial knot invariants \cite{stoimenow}, linking forms \cite{Lickorish} and $3$-manifold theory \cite{gordon-luecke}.    
In this article we propose to add yet another tool to this list by using the rational Witt class $\varphi (K)$ to extract information about $u(K)$. The rational Witt class $\varphi (K)$ is associated to an oriented knot $K\subset S^3$ and takes values in the Witt ring $W(\qq) \cong \mathbb Z \oplus \mathbb Z_2^\infty \oplus \mathbb Z_4^\infty$ of the field $\qq$ of rational numbers.\footnote{Here and below, we write $\mathbb Z_p$ to mean $\mathbb Z/p\mathbb Z$ while we use $\mathbb Z_p^\infty $ as a shorthand for $\oplus _{i=1}^\infty \mathbb Z_p$.} As a commutative ring,  $W(\qq)$ is obtained by applying the Grothendieck group construction to the Abelian semiring 
of isomorphism classes of non-degenerate, symmetric, bilinear forms on finite dimensional rational vector spaces. The operations on the latter are given by direct sums and tensor products of vector spaces along with summing and multiplying their bilinear forms. The Witt ring $W(\qq)$ is well understood and we describe it in some detail in Section \ref{section-on-witt-rings}. For the time being, we content ourselves with saying that $W(\qq)$ is generated by $1$-dimensional forms $\langle a \rangle$, $a\in \dot{\qq}$ where $\langle a\rangle :\mathbb Q \times \qq \to \mathbb Q$ is the unique bilinear form that sends $(1,1)$ to $a$ (and, as usual, $\dot\qq = \qq-\{0\}$). Thus, given a non-degenerate, symmetric, bilinear form $q:\mathbb Q^n \times \qq^n \to \qq$, there exist rational numbers $a_1,...,a_n \in \dot\qq$ such that $q=\langle a_1\rangle \oplus ... \oplus \langle a_n\rangle$, i.e. such that $q((x_1,...,x_n),(y_1,...,y_n)) = a_1x_1 y_1+...+a_nx_n y_n$. 

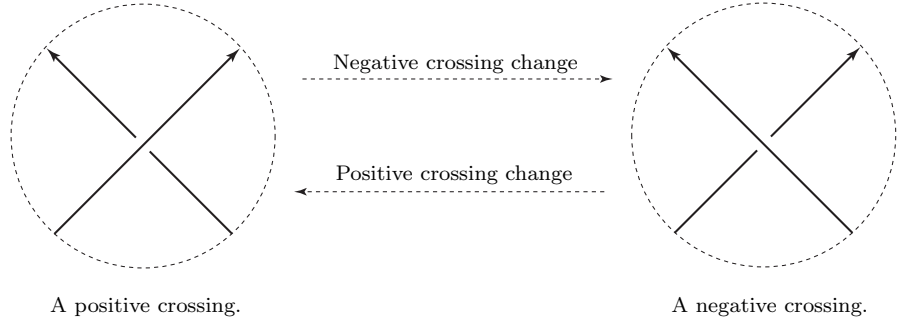
\begin{figure}[htb!] 
\centering
\includegraphics[width=12cm]{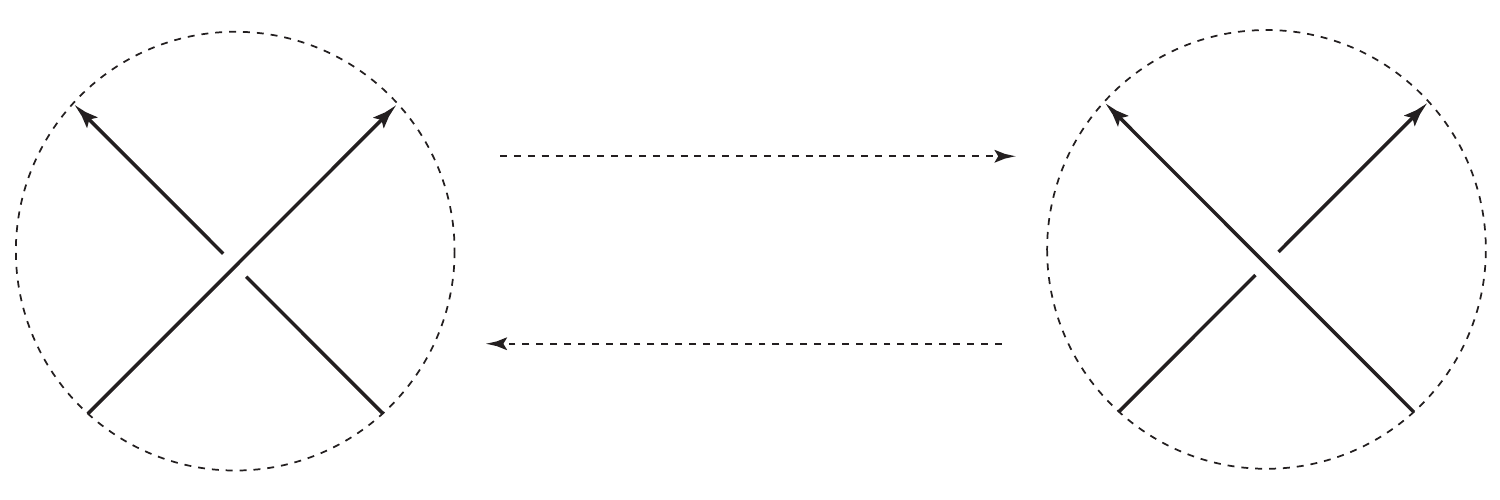}
\put(-216,82){\tiny Negative crossing change}
\put(-215,40){\tiny Positive crossing change}
\put(-323,-10){\tiny A positive crossing.}
\put(-88,-10){\tiny A negative crossing.}
\caption{Our convention for positive and negative crossings as well as positive and negative crossing changes. }  \label{pic1}
\end{figure}

Given an oriented knot $K$ in $S^3$, we shall label crossings in a projection of $K$ as {\em positive} or {\em negative} according to the usual convention, see Figure \ref{pic1}. Crossing changes themselves shall be similarly labeled as {\em positive} or {\em negative} according to whether they change a negative crossing to a positive one or vice versa, see again Figure \ref{pic1}. 

Our main results, Theorem \ref{main1} and its various corollaries, are founded on the observation that $\varphi (K)$ changes rather predictably when $K$ undergoes a single crossing change. This phenomenon is described in the next statement. 
\begin{theorem} \label{main1}
Let $K_+$ be a knot obtained from the knot $K_-$ by a positive crossing change. Then the rational Witt classes of $K_+$ and $K_-$ are related as follows, depending on how their signatures $\sigma (K_\pm)$ compare: 
$$
\begin{array}{rl}
\varphi(K_+)&=
\left\{
\begin{array}{ll}
\varphi(K_-) \oplus \left\langle\textstyle \frac{2\det K_+}{\det K_-} \right\rangle \oplus \langle -2\rangle & \quad ; \quad \sigma (K_+) = \sigma (K_-)  \cr  & \cr
\varphi(K_-) \oplus \left\langle \textstyle -\frac{2\det K_+}{\det K_-} \right\rangle \oplus \langle -2\rangle & \quad ; \quad \sigma (K_+) = \sigma (K_-) -2 
\end{array}\right.  \cr & \cr
\varphi(K_-)&=
\left\{
\begin{array}{ll}
\varphi(K_+) \oplus \left\langle\textstyle -\frac{2\det K_-}{\det K_+} \right\rangle \oplus \langle 2\rangle & \quad ; \quad \sigma (K_-) = \sigma (K_+)  \cr  & \cr
\varphi(K_+) \oplus \left\langle \textstyle \frac{2\det K_-}{\det K_+} \right\rangle \oplus \langle 2\rangle & \quad ; \quad \sigma (K_-) = \sigma (K_+) +2 
\end{array}\right.
\end{array}
$$
The conditions on $\sigma (K_\pm)$ stated on the right-hand sides above, cover all possible cases. 
\end{theorem}
A similar, and rather beautiful formula for how the algebraic concordance class of a knot changes under a crossing switch, was found by S.-G. Kim and C. Livingston in \cite{livingston-kim}.  

As the rational Witt class of the unknot is trivial, the Theorem \ref{main1} gives restrictions on what $\varphi (K)$ can be if $K$ has a given unknotting number. While such restriction exist regardless of the value of $u(K)$, they are easiest to state, and have proven most effective, when $u(K) = 1$. 
\begin{corollary} \label{coro1}
Let $K$ be a knot with unknotting number $1$. Then the rational Witt class $\varphi (K)$ of $K$ must be as follows, depending on whether $K$ can be unknotted by a positive or a negative crossing change.
\begin{itemize}
\item[a)] If $K$ can be unknotted by a positive crossing change, then 
$$\phantom{iiiii} \varphi (K) = \left\{
\begin{array}{rl} 
\displaystyle \langle 2\det K \rangle \oplus \langle 2 \rangle  & \quad ; \quad \sigma (K) = 2 \cr & \cr 
\displaystyle \langle -2\det K \rangle \oplus \langle 2 \rangle & \quad ; \quad \sigma (K) = 0 
\end{array}
\right. 
$$
\item[b)] If $K$ can be unknotted with a negative crossing change, then 
$$\varphi (K) = \left\{
\begin{array}{rl} 
\displaystyle \langle 2\det K \rangle \oplus \langle -2 \rangle & \quad ; \quad \sigma (K) = 0 \cr & \cr 
\displaystyle \langle -2\det K \rangle \oplus \langle -2 \rangle & \quad ; \quad \sigma (K) = -2  
\end{array}
\right. 
$$
\end{itemize}
As before, $\sigma (K)$ denotes the signature of $K$. 
\end{corollary}
The next corollary provides similar constraints on the rational Witt class $\varphi (K)$ of a knot $K$ with $u(K)=2$. It is inherently weaker than Corollary \ref{coro1} in that it involves information about the knot $L$ obtained from $K$ after only one crossing change, a knot which one generally knows little about.  

\begin{corollary} \label{coro2}
Let $K$ be a knot with unknotting number $2$ and let $L$ be the knot obtained from $K$ after only a single crossing change. Then the rational Witt class $\varphi (K)$ is determined by $\det K$, $\det L$ and $\sigma (K)$ and the type of crossing changes involved, as indicated below. To reduce the number of cases to state, we make the assumption that $\sigma (K) \le 0$.\footnote{The assumption of $\sigma (K) \le 0$ in Corollary \ref{coro2} can always be achieved by, if necessary, replacing $K$ by its mirror image $\bar K$. Clearly $u(\bar K) = u(K)$ while $\sigma (\bar K) = -\sigma (K)$ and $\varphi (\bar K) = -\varphi (K)$.} 
\begin{itemize}
\item[a)] If $K$ can be unknotted by two negative crossing changes, then 
$$\varphi (K) = \left\{ 
\begin{array}{rl}
\langle -2 \det K \det L\rangle \oplus \langle -2 \det L\rangle \oplus  \langle -1\rangle \oplus \langle -1\rangle  & ; \sigma (K) = -4 \cr
& \cr
\langle \pm 2 \det K \det L\rangle \oplus \langle \mp 2 \det L\rangle \oplus  \langle -1\rangle \oplus \langle -1\rangle & ; \sigma (K) = -2 \cr
 & \cr 
\langle 2 \det K \det L\rangle \oplus  \langle 2 \det L\rangle \oplus  \langle -1\rangle \oplus \langle -1\rangle & ; \sigma (K) = 0  
\end{array}
\right.
$$
The signs in the second line have to be chosen consistently either as $(+, -)$ (if $\sigma (L) = -2$) or as $(-, +)$ (if $\sigma (L) = 0$). 
\item[b)] If $K$ can be unknotted by one positive and one negative crossing change, then 
$$\varphi (K) = \left\{ 
\begin{array}{cl}
\langle -2 \det K \det L\rangle \oplus  \langle -2 \det L\rangle& ; \sigma (K) = -2 \cr
 & \cr 
\langle \pm 2 \det K \det L\rangle  \oplus   \langle \mp 2 \det L\rangle & ; \sigma (K) = 0  
\end{array}
\right.
$$
Here too the signs in the $\sigma (K)=0$ case have to be chosen consistently. The choice of $(+, -)$ corresponds to the case where $L$ is either obtained from $K$ by a negative crossing change and $\sigma (L) = \sigma (K) + 2$ or $L$ is obtained from $K$ by a positive crossing change and $\sigma (L) = \sigma (K)$. The choice $(-,+)$ represents the other two possibilities.  
\item[c)] If $K$ can be unknotted with two positive crossing changes then $\sigma (K) = 0$ and 
$$ \textstyle \varphi (K) = \langle - 2 \det K \det L\rangle \oplus  \langle -2 \det L\rangle \oplus  \langle 1\rangle \oplus \langle 1\rangle $$
\end{itemize} 
%
\end{corollary}
Our techniques apply equally well to knots with higher unknotting numbers. However, the indeterminacy of $\varphi (K)$ of a knot $K$ with $u(K) = n$ grows with $n$ in that it involves the determinants of all the knots that $K$ \lq\lq goes through on its way to the unknot\rq\rq. This phenomenon substantially diminishes the usefulness  of our approach, more so since the number of cases describing $\varphi (K)$ grows with $n$ as well. We list the next corollary more as an illustration of our methods rather than a tool we deem practically useful.   
\begin{corollary} \label{coro3}
Let $K$ be a knot with signature $-2n$ and with $u(K) =n\ge 1$. Let $L_i$ be the knot obtained from $K$ by changing $i-1$ of the $n$ crossings (e.g. $L_1$ is just $K$ while $L_{n+1}$ is the unknot). Then 
$$ \varphi (K) = \bigoplus _{i=1}^n \left(   \langle -2 \det L_{i+1}\, \det L_i \rangle \oplus \langle -2\rangle \right)
$$
\end{corollary}
We remark that both the $\sigma (K) = -2$ case in Corollary \ref{coro1} and the $\sigma (K) = -4$ case in Corollary \ref{coro2} follow from Corollary \ref{coro3} after recognizing that the equality
$$
\bigoplus _{i=1}^n \langle -2\rangle= 
\left\{ 
\begin{array}{cl}
\langle-2\rangle \oplus \left(\bigoplus _{i=1}^{n-1} \langle -1\rangle \right)   & \quad ; \quad n \text{ is odd} \cr 
& \cr
\bigoplus _{i=1}^{n} \langle -1\rangle  & \quad ; \quad n \text{ is even} 
\end{array}
\right.
$$
holds in $W(\qq)$ for each $n\in \mathbb N$. 
\subsection{Applications and examples}
The main utility of Theorem \ref{main1} and its corollaries is to provide obstructions for a given knot $K$ to satisfy the equation $u(K) = n$. To a large degree, our emphasis shall be on the case $n=1$. 

We start this section by subjecting $3$ different families of knots to Corollary \ref{coro1}. The first two of these families are finite and consist of $151$ knots with $11$ crossings and $100$ knots with $12$ crossings respectively. The third family is that of $n$-stranded pretzel knots for various $n\ge 3$. We then apply Corollary \ref{coro2} to the knot $K=10_{47}$ which, at the time of this writing, has unknown unknotting number \cite{knotinfo2} (though it is either $2$ or $3$). 

As Witt rings live at the interface of number theory, algebra and -- to a minor degree -- topology, the  reader will likely detect a number theoretic flair in many of our subsequent statements, especially those regarding pretzel knots. 
\subsubsection{Eleven crossing knots} 
We consider the family of alternating and non-alternating $11$ crossing knots  $11a_x$ and $11n_y$ with $x$ and $y$ ranging through the following parameter sets, organized by signature (see Remark \ref{remark-color-highlighting} below for an explanation of the color highlighting).
\begin{equation}  \label{the-eleven-crossing-knots}
\begin{array}{l}
x\in 
\left\{
\begin{array}{ll}
\left\{ 
\begin{minipage}[c]{10.3cm} 
\colorbox{green}{7}, \colorbox{green}{33}, 51, 55, 92, 108, 131, \colorbox{green}{137}, 155, 158, 162, 196, 199, 217, 218, 
\colorbox{green}{219}, 221, 248, 273, \colorbox{green}{296}, \colorbox{green}{297}, 301, 305, 312, 322, 324, 325, 331 
\end{minipage} 
\right\} & ;  \sigma (11a_x) = 2  \cr
& \cr 
\left\{ 
\begin{minipage}[c]{10.3cm} 
4, 5, \colorbox{green}{16}, 36, \colorbox{yellow}{37}, 39, 58, 87, 103, 109, 112, \colorbox{yellow}{128}, 135, 153, 164, 
165, 169, \colorbox{green}{170}, 201,  \colorbox{yellow}{214}, \colorbox{yellow}{228}, 249, 270, \colorbox{green}{274}, \colorbox{yellow}{278}, \colorbox{yellow}{285}, \colorbox{green}{288},  
303, \colorbox{yellow}{313}, 315, 317, \colorbox{yellow}{332}, 350
\end{minipage} 
\right\} & ;  \sigma (11a_x) = 0  \cr
& \cr
\left\{ 
\begin{minipage}[c]{10.3cm}  
1, 6, 21, 23, 32, 42, \colorbox{green}{45}, 46, 50, 61, 97, \colorbox{green}{99}, \colorbox{green}{107}, 118, 125, 133, 
134, \colorbox{green}{148}, \colorbox{green}{163}, 171, 172, 181, 197, \colorbox{green}{202}, \colorbox{green}{239}, 258, 268, 269, 271, 277, 
279, \colorbox{green}{281}, 284, 286, \colorbox{green}{314}, 327, 349, 352, 362 
\end{minipage} 
\right\} & ;  \sigma (11a_x) = -2
\end{array} \right. 
\cr \cr \cr
y \in 
\left\{
\begin{array}{ll}
\left\{ 
\begin{minipage}[c]{10.3cm}   
3, 17, \colorbox{green}{58}, 91, 92, 102, \colorbox{green}{113}, 122, 127, 129, \colorbox{green}{140}, 170  
\end{minipage} 
\right\} & ;  \sigma (11n_y) = 2  \cr
& \cr
\left\{ 
\begin{minipage}[c]{10.3cm}   
49, 51, 83, \colorbox{yellow}{94}, \colorbox{yellow}{115}, 116, \colorbox{yellow}{119}, 132, 139, \colorbox{yellow}{141}, 142, 157, \colorbox{green}{165}, 172, \colorbox{yellow}{179}, \colorbox{yellow}{182}
\end{minipage} 
\right\} & ;  \sigma (11n_y) = 0  \cr
& \cr
\left\{ 
\begin{minipage}[c]{10.3cm}   
\colorbox{green}{15}, \colorbox{green}{29}, 54, 60, \colorbox{green}{79}, 112, \colorbox{green}{117}, 120, 128, 138, 146, 148, 150, \colorbox{green}{155}, 160, 161, 162, 
\colorbox{green}{163}, 166, 167, 168, 177, 178
\end{minipage} 
\right\} & ;  \sigma (11n_y) = -2
\end{array} \right. 
\end{array}
\end{equation}
Prior to the recent results by J. Greene \cite{Greene}, these were the knots with $11$ crossings whose unknotting numbers were unknown but were either $1$ or $2$, according to KnotInfo \cite{knotinfo1}. The exception to this were the knots 
$$11a_{45}, 11a_{137}, 11a_{197}, 11a_{202}, 11a_{362}\quad \quad \text{ and } \quad \quad 11n_{141}, 11n_{148}$$
which were listed as having unknotting number either $1$, $2$ or $3$. 

While Greene's results \cite{Greene} show that none of these knots can have unknotting number $1$, and his results thereby subsume our findings, we nevertheless list here the outcome of applying Corollary \ref{coro1}
to the above knots as an illustration of its efficacy and in the hopes that the reader may appreciate an alternate and independent proof of some of the results from \cite{Greene}.
\begin{corollary} \label{theorem-results-for-11-crossing-knots}
Consider the knots $11a_x$ and $11n_y$ with $x$ and $y$ as in \eqref{the-eleven-crossing-knots}. 
\begin{itemize}
\item[(a)] Each of the knots $11a_x$ and $11n_y$, with $x$ and $y$, from 
$$
\begin{array}{l}
x\in 
\left\{
\begin{array}{ll}
\left\{ 
\begin{minipage}[c]{6.9cm} 
7, 33, 137, 219, 296, 297
\end{minipage} 
\right\} & \quad  ;  \quad   \sigma (11a_x) = 2  \cr
& \cr 
\left\{ 
\begin{minipage}[c]{6.9cm} 
16, 170, 274, 288
\end{minipage} 
\right\} &  \quad  ;  \quad  \sigma (11a_x) = 0  \cr
& \cr
\left\{ 
\begin{minipage}[c]{6.9cm}  
45, 99, 107, 148, 163, 202, 239, 281, 314
\end{minipage} 
\right\} &  \quad  ;  \quad   \sigma (11a_x) = -2
\end{array} \right. 
\cr \cr \cr
y \in 
\left\{
\begin{array}{ll}
\left\{ 
\begin{minipage}[c]{6.9cm}   
58, 113, 140
\end{minipage} 
\right\} &  \quad  ;  \quad   \sigma (11n_y) = 2  \cr
& \cr
\left\{ 
\begin{minipage}[c]{6.9cm}   
165
\end{minipage} 
\right\} &  \quad  ;  \quad  \sigma (11n_y) = 0  \cr
& \cr
\left\{ 
\begin{minipage}[c]{6.9cm}   
15, 29, 79, 117, 155, 163 
\end{minipage} 
\right\} &  \quad  ;  \quad  \sigma (11n_y) = -2
\end{array} \right. 
\end{array}
$$ 
have unknotting number 2, with the possible exception of $11a_{45}$, $11a_{137}$ and $11a_{202}$ which have unknotting number at least $2$. 
\item[(b)] None of the signature zero knots $11a_x$ or $11n_y$ with $x$ and $y$ from 
$$x\in \{37, 214, 278, 313    \} \quad \text{ and } \quad y\in \{179 \}$$
can be unknotted with a single negative crossing change. 
\item[(c)] None of the signature zero knots $11a_x$ or $11n_y$ with $x$ and $y$ from
$$x\in \{128, 228, 285, 332  \} \quad \text{ and } \quad y\in \{94, 115, 119, 141, 182 \}$$
can be unknotted with a single positive crossing change. 
\end{itemize}
\end{corollary}
\begin{remark} \label{remark-color-highlighting}
The knots from part (a) of the preceding corollary have been shaded green in \eqref{the-eleven-crossing-knots} while those from parts (b) and (c) are represented by a yellow shading.   We note that of the 151 knots from \eqref{the-eleven-crossing-knots}, Corollary \ref{theorem-results-for-11-crossing-knots} provides unknotting information for $43$ of them. 
\end{remark}
\subsubsection{Twelve crossing knots.}
At the time of this writing, no data concerning the unknotting numbers of $12$ crossings knots is available on KnotInfo \cite{knotinfo2}. In order to start collecting such data, we have applied Corollary \ref{coro1} to the first 50 alternating and the first 50 non-alternating 12 crossing knots. Considering that there are 1288 alternating and 888 non-alternating 12 crossing knots total, this is but a very modest beginning to a somewhat daunting  program. 

Among the knots $12a_x$ and $12n_y$ with $1\le x,y \le 50$, the following have signature greater than $2$ or less than $-2$  
\begin{equation} \label{knots-to-be-excluded-from-12-crossing-theorem}
\begin{array}{rl}
x\in & \{ 11,  21, 24, 26, 34, 35, 36, 37,  46, 50\}\cr \cr
y\in & \{ 6, 8, 16, 37 \} 
\end{array}
\end{equation}
Since such knots have unknotting number at least $2$ (by virtue of \eqref{lower-bounds-on-uK}), we exclude these from the next corollary. 
\begin{corollary} \label{theorem-results-for12-crossing-knots}
Consider the knots $12a_x$ and $12n_y$ with $1\le x, y  \le 50$ but with the exception of those $x$ and $y$ listed in \eqref{knots-to-be-excluded-from-12-crossing-theorem}. 
\begin{itemize}
\item[(a)] The knots $12a_x$ and $12n_y$ with   
$$
\begin{array}{l}
x\in 
\left\{
\begin{array}{ll}
\left\{ 
\begin{minipage}[c]{1.9cm} 
29, 32, 39
\end{minipage} 
\right\} & \quad  ;  \quad   \sigma (12a_x) = 2  \cr
& \cr 
\left\{ 
\begin{minipage}[c]{1.9cm} 
48
\end{minipage} 
\right\} &  \quad  ;  \quad  \sigma (12a_x) = 0  \cr
& \cr
\left\{ 
\begin{minipage}[c]{1.9cm}  
9 
\end{minipage} 
\right\} &  \quad  ;  \quad   \sigma (12a_x) = -2
\end{array} \right. 
\cr \cr \cr
y \in 
\left\{
\begin{array}{ll}
\left\{ 
\begin{minipage}[c]{1.9cm}   
10, 15,  17
\end{minipage} 
\right\} &  \quad  ;  \quad   \sigma (12n_y) = 2  \cr
& \cr
\left\{ 
\begin{minipage}[c]{1.9cm}   
27, 33
\end{minipage} 
\right\} &  \quad  ;  \quad  \sigma (12n_y) = -2
\end{array} \right. 
\end{array}
$$ 
have unknotting number at least $2$. 
\item[(b)] None of the signature zero knots $12a_x$ and $12n_y$ with 
$$ x\in \{16\} \quad \quad \text{ and } \quad \quad y\in \{1, 22, 30\}$$
can be unknotted with a single negative crossing change. 
\item[(c)] None of the signature zero knots $12a_x$ and $12n_y$ with 
$$ x\in \{1, 13, 15, 23, 30, 33, 43\} \quad \quad \text{ and } \quad \quad y\in \{28, 34, 35, 39\}$$
can be unknotted with a single positive crossing change.
\end{itemize}
\end{corollary}
Corollary \ref{theorem-results-for12-crossing-knots} produces unknotting information for $25$ of the $86$ examined knots.  
\subsubsection{Pretzel knots} \label{introductory-section-on-pretzel-knots}
\begin{figure}[htb!] 
\centering
\includegraphics[width=5cm]{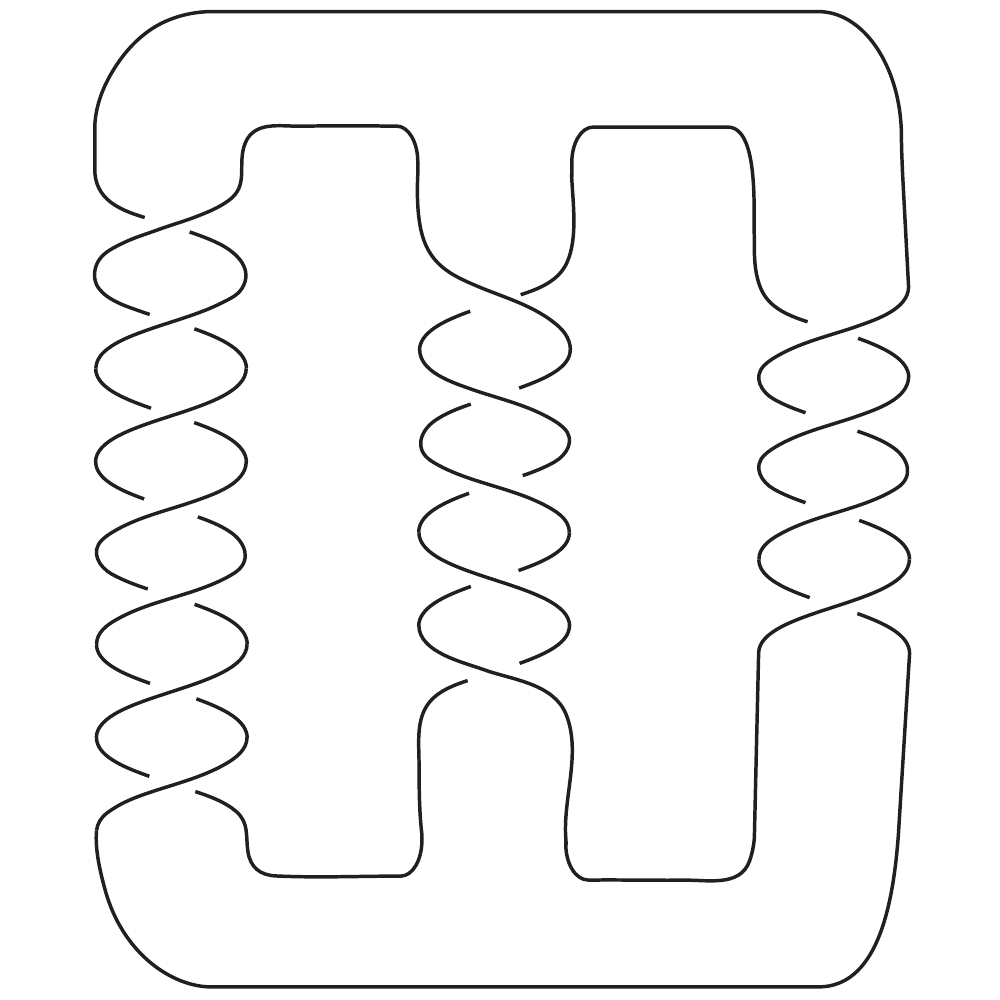}
\caption{The pretzel knot $P(7,-5,4)$.}  \label{pic4}
\end{figure}

For an integer $n\ge 3$ and for nonzero integers $p_1, ..., p_n$, we let $P(p_1,...,p_n)$ denote the corresponding $n$-stranded pretzel knot/link. It is obtained by taking $n$ pairs of parallel strands, introducing $|p_i|$ half-twists into the $i$-th pair (with $p_i>0$ giving right-handed and $p_i<0$ left-handed half-twists), and closing up the strands with $n$ pairs of bridges. Figure \ref{pic4} shows the example $P(7, -5, 4)$.  
In order for $P(p_1,...,p_n)$ to be a knot, at most one of $p_1,...,p_n$ can be even. In fact, if $n$ itself is even, then precisely one of $p_1, ..., p_n$ needs to be even. We shall assume these parity conditions to be satisfied throughout. 

The computation of $\varphi (P(p_1,...,p_n))$ in \cite{jabuka1}, for any choice of parameters $p_1,...,p_n$, provides a fertile testing ground for Theorem \ref{main1} and Corollary \ref{coro1}. We list here only a few select applications and examples, leaving a more comprehensive exploration of unknotting numbers of pretzel knots, for a future occasion. We start with the following remark. 
\begin{remark} \label{about-unknotting-of-two-bridge-knots}
In \cite{kawauchi}, A. Kawauchi showed that a pretzel knot $P(p_1,...,p_n)$ is a two-bridge knot precisely when at most two of $p_1,...,p_n$ differ from $\pm 1$.  As T. Kanenobu and H. Murakami \cite{kanenobu-murakami}  determined all two-bridge knots with unknotting number $1$, we omit such knots from our applications below. 

In \cite{kobayashi}, T. Kobayashi showed that the pretzel knots $P(p_1,p_2,p_3)$ with $p_1,p_2,p_3$ odd and with $u(P(p_1,p_2,p_3))= 1$, are precisely those non-trivial knots for which $\{a,b\}\subset \{p_1,p_2,p_3\}$ where $\{a,b\}$ is either $\{\pm 1, \pm 1\}$ or $\{\pm 3, \mp 1\}$. We shall therefore also exclude such knots from our examples. 
\end{remark}
As Corollary \ref{coro1} is sensitive to signatures, we remark that the signatures of pretzel knots $P(p_1,...,p_n)$ have also been computed in full generality in \cite{jabuka1}.  
\begin{corollary} \label{pretzel1}
Consider the pretzel knot $P(p_1,p_2,p_3)$ with $p_1, p_2$ odd and with $p_3$ even. 
If $p_1\ge 7$, $p_3 > -\frac{p_1(4-p_1)}{4}$ and the equality 
\begin{align} \nonumber
\langle -1\rangle \oplus \langle -2\rangle  \oplus \langle 4p_3 - p_1(p_1-4)\rangle 
=\langle -2(4p_3 - p_1(p_1-4))\rangle 
\end{align} 
fails to hold in $W(\qq)$, then $u(P(p_1,4-p_1,p_3))\ge 2$. 
\end{corollary}
There are many examples meeting the hypothesis of Corollary \ref{pretzel2}. Here are a couple. 
\begin{example} \label{example-one-for-pretzels-with-p3-even}
Consider the pretzel knot $P(7,-3,p_3)$ with $p_3\ge 6$ an even integer. If there exists a prime $p$ dividing $4p_3-21$ with an odd power and such that $-2$ is not a square in $\mathbb Z_p$, then $u(P(7,-3,p_3))\ge 2$. For example, any of $p_3=2k\cdot 7^{\ell+1}$ with $k,\ell \in \mathbb N$, satisfies these conditions (with $p=7$).
\end{example}
\begin{example} \label{example-two-for-pretzels-with-p3-even}
The unknotting number of $P(17,-13,p_3)$ with $p_3=15+(2k+1)\cdot 23^{\ell+1}$ and $k,\ell \in \mathbb N$, is at least $2$. 
\end{example}
\begin{corollary} \label{pretzel2}
Let $p>0$ be an odd integer and consider the $4$-stranded pretzel knot $P(p,p,p,-3p-1)$. If the  equality 
$$
\langle 1 \rangle \oplus  \langle p\rangle \oplus  \langle p\rangle \oplus  \langle p\rangle \oplus \langle -3p-1\rangle \oplus \left\langle \textstyle - \frac{8p+3}{p(3p+1)}  \right\rangle =  \langle 2\rangle \oplus \langle 2(8p+3) \rangle 
$$
fails to hold in $W(\qq)$, then $u(P(p,p,p,-3p-1))\ge 2$. 
\end{corollary}
Many choices of $p$ are possible in Corollary \ref{coro2}. Here is an infinite family of such choices. 
\begin{example} \label{example-for-4-stranded-pretzels}
Taking $p=2+(2k+1)\cdot 19^{\ell+1}$ with $k,\ell \in \mathbb N$, meets the conditions of Corollary \ref{pretzel2}. Consequently, each of the corresponding knots $P(p,p,p,-3p-1))$ has unknotting number at least $2$. 
\end{example}
Given a pretzel knot $P(p_1,...,p_n)$ and an odd integer $p$, we say that the knot $P(p_1,...,p_m,p,p_{m+1},...,p_\ell,-p,p_{\ell+1},...,p_n)$ was obtained from $P(p_1,...,p_n)$ by {\em upward stabilization} (a term already introduced in \cite{jabuka1}). With this in mind, we have:
\begin{corollary} \label{pretzel3}
If $K=P(p_1,...,p_n)$ is a pretzel knot for which the equalities for $\varphi (K)$ from Corollary \ref{coro1} fail (so that $u(K)\ge 2$), then the same is true for any pretzel knot $L$ obtained from $K$ by a finite number of upward stabilizations.  Consequently, for any such $L$ one has $u(L)\ge 2$.  
\end{corollary}
Combining this corollary with previous examples, supplies pretzel knots with an arbitrarily high number of strands and with unknotting number at least $2$. For instance, 
$$u(P(7,-3,14, p_1,-p_1,p_2,-p_2,...,p_m,-p_m))\ge 2$$
for any choice of odd integers $p_1,...,p_m$.  
\subsubsection{Obstructing unknotting number $2$} 
As already mentioned, Corollary \ref{coro2} is inherently weaker than Corollary \ref{coro1} as it involves the unknown quantity $\det L$. Even so, it is still possible to gain some unknotting information from it. To demonstrate this, we consider the knot $K=10_{47}$ which has signature $4$, determinant $41$ and, as of this writing, has unknotting number $2\le u(10_{47}) \le 3$, according to KnotInfo \cite{knotinfo2}. Applying Corollary \ref{coro2} to 
this knot, we find:
\begin{corollary} \label{corollary-for-unknotting-number-two}
Suppose $K=10_{47}$ can be unknotted by $2$ crossing changes and let $L$ be the knot obtained from $K$ by a single of these crossing changes. If $L$ has $9$ or fewer crossings, then $L$ must be contained in the list of $12$  knots (out of $84$ knots with $9$ or fewer crossings, not counting mirror images):
$$\bar 3_1, \bar 5_2, \bar 6_2, \bar 7_2, \bar 7_6, \bar 8_{11}, \bar 8_{21}, \bar 9_2, \bar 9_{12}, 9_{26}, 9_{39}, 9_{42}. 
$$
A bar on top of a knot indicates its mirror image.
\end{corollary}
\vskip3mm
We finish this section by pointing out that all of our applications of Theorem \ref{main1} made the choice of either $K_+ =$ unknot or $K_-=$ unknot. The usefulness of Theorem \ref{main1} certainly stretches beyond this. We leave it as an exercise for the motivated reader to verify, for example, that the knots $8_4$ and $9_{19}$ cannot be gotten from one another by a single crossing change. 
\subsection{Organization}
The remainder of this article is organized into 6 sections. Section \ref{section-on-witt-rings} provides background material on Witt rings with a special emphasis on the Witt ring of the rationals. Section \ref{section-phi-of-k-under-a-crossing-change}  defines the rational Witt class $\varphi (K)$ associated to a knot $K$ and explores how the former changes when $K$ is altered by a single crossing change. Doing so enables us to prove Theorem \ref{main1}. Section \ref{section-the-proofs-of-the-theorems} supplies the proofs for Corollaries \ref{coro1}, \ref{coro2} and \ref{coro3} while Section \ref{section-about-the-low-crossing-knots} explains how the results from Corollaries \ref{theorem-results-for-11-crossing-knots}, \ref{theorem-results-for12-crossing-knots} and \ref{corollary-for-unknotting-number-two} were obtained. Section \ref{section-on-pretzel-knots} provides proofs of our claims concerning pretzel knots while the final Section \ref{section-work-of-lickorish} provides a comparison of our work to that of R. Lickorish from \cite{Lickorish}.

\vskip3mm
{\bf Acknowledgements } During the preparation of this article, I have enjoyed and benefited from conversations with Brendan Owens, whose input I gratefully acknowledge. I would also like to thank Josh Greene for generously sharing his results from \cite{Greene} and for providing helpful comments on an earlier version of this work. Additional thanks are due to Chuck Livingston and Swatee Naik for many stimulating conversations about Witt rings. 
\section{Background material on Witt rings} \label{section-on-witt-rings}
This section reviews some of the basic algebra underlying the definition of Witt rings $W(\mathbb F)$ over arbitrary fields $\mathbb F$. We then focus in on the case of $\mathbb F = \qq$ and give a completely explicit description of the isomorphism $W(\qq)\cong \mathbb Z \oplus \mathbb Z_2^\infty \oplus \mathbb Z_4^\infty$ which was already mentioned in the introduction. For more information we advise the interested reader to consider the sources \cite{gerstein, husemoller-milnor, lam, omeara}. 

To begin with, let us fix a field $\mathbb F$ and let $\mathcal B_\mathbb F$ be the set of isomorphism classes of symmetric, bilinear, non-degenerate forms over finite dimensional $\mathbb F$-vector spaces. Thus, an element of $\mathcal B_\mathbb F$ is a pair $(V,B)$ where $V$ is a finite dimensional $\mathbb F$-vector space and $B:V\times V\to \mathbb F$ is a symmetric, bilinear and non-degenerate form where by the latter we mean that if $B(x,y) = 0$ for all $y\in V$,  then $x=0$. The set $\mathcal B_\mathbb F$ becomes an Abelian semiring (often, somewhat humorously, referred to as an {\em Abelian rig}) under the operations $\oplus$ and $\otimes$ given by  
$$(V_1, B_1) \oplus (V_2, B_2) = (V_1\oplus V_2, B_1+B_2) \quad \text{ \& } \quad (V_1, B_1) \otimes (V_2, B_2) = (V_1\otimes _\mathbb F V_2, B_1\cdot B_2)$$
When the importance of $V$ is minor and the danger of confusion little, we will only write $B$ to mean $(V,B)$ and likewise $B_1\oplus B_2$ to mean $(V_1,B_1) \oplus (V_2, B_2)$. 
 
Recall that {\em Grothendieck's group construction} turns an Abelian semigroup $(G,+)$ into an Abelian group by considering the set $(G\times G)/\hspace{-2mm}\sim$ where $\sim$ is the equivalence relation defined by 
$$(x_1,y_1) \sim (x_2,y_2) \quad \text{ if } \quad x_1+y_2 = x_2+y_1$$
(intuitively we should regard $(x,y)$ as representing $x-y$, even though the latter is of course not defined). With respect to the addition $(x_1,y_1)+(x_2,y_2) = (x_1+x_2, y_1+y_2)$ on $(G\times G)/\hspace{-2mm}\sim$, the inverse of $(x,y)$ is then given by $(y,x)$. The semigroup $G$ injects naturally into $(G\times G)/\hspace{-2mm}\sim$ by sending $x$ to $(x,0)$. If $G$ has the structure of an Abelian semiring, $(G\times G)/\hspace{-2mm}\sim$ itself becomes an Abelian ring. 

With this understood, here is the definition of the Witt ring. 
\begin{definition} 
The Witt ring $W(\mathbb F)$ associated to the field $\mathbb F$, is the Abelian ring obtained by applying the Grothendieck construction to the Abelian semiring $(\mathcal B_\mathbb F, \oplus, \otimes)$. 
\end{definition}
As is customary, we shall use $\dotf$ to denote $\mathbb F-\{0\}$. Given an element $a\in  \dotf$, let $\langle a \rangle $ denote the unique bilinear, symmetric, non-degenerate form on $\mathbb F \times \mathbb F$ which sends $(1,1)$ to $a\in \mathbb F$. Note that 
$$\langle a \rangle = \langle a\cdot d^2\rangle$$
for any choice of $d\in \dotf$ since $f:(\mathbb F, \langle a\cdot d^2\rangle) \to (\mathbb F, \langle a\rangle)$ given by $f(\lambda ) = d\cdot \lambda$ is an isomorphism of bilinear forms. We will often tacitly rely on the equality $\langle a \rangle = \langle a\cdot d^2\rangle$ in the remainder of the article. 

The next theorem is basic and can be found in each of \cite{gerstein, husemoller-milnor, lam, omeara}. 
\begin{theorem} \label{presentation-of-the-rational-witt-ring}
For any field $\mathbb F$, the Witt ring $W(\mathbb F)$ is generated by the set $\{ \langle a\rangle \, |\, a\in \dotf \}$. A presentation of $W(\mathbb F)$ as a commutative ring is obtained by adding the next relators to these generators:
\begin{align} \nonumber
(R1) \quad \quad & \langle 1 \rangle \oplus \langle -1\rangle \cr
(R2) \quad \quad & \langle a \rangle \otimes \langle b \rangle \oplus \langle -a\cdot b\rangle & \quad \quad \quad a,b \in \mathbb \dotf \cr 
(R3) \quad \quad & \langle a+b \rangle \oplus \langle ab (a+b)\rangle \oplus \langle -a\rangle \oplus \langle -b \rangle & \quad \quad \quad a,b \in \mathbb \dotf
\end{align}
\end{theorem}
For the next discussion, we assume that $\text{char } \mathbb F \ne 2$. 
The {\em hyperbolic form} over the field $\mathbb F$ is the $2$-dimensional bilinear form $(\mathbb F^2, H)$ where $H$, with respect to the standard basis $\{ e_1=(1,0), e_2=(0,1) \}$ of $\mathbb F^2$, is represented by the matrix 
$$H = \left[
\begin{array}{cc}
0 & 1 \cr
1 & 0 
\end{array}
\right]
$$
With respect to the basis $\{f_1, f_2\}$ of $\mathbb F^2$, given by $f_1=\frac{1}{2} e_1+e_2$, $f_2=-\frac{1}{2} e_1 +e_2 $, the hyperbolic form $H$ is represented by the matrix
$$  H = \left[
\begin{array}{cr}
1 & 0 \cr
0 & -1 
\end{array}
\right] = \langle 1 \rangle \oplus \langle -1 \rangle =0 \in W(\mathbb F)
$$
showing that it equals zero in $W(\mathbb F)$. 

Hyperbolic forms are rather special, indeed, they arise as summands of all \lq\lq isotropic forms\rq\rq. A form $(V,B)$ is called {\em isotropic} if there exists a non-zero vector $v\in V$ with $B(v,v)=0$, otherwise $(V,B)$ is called {\em anisotropic}. Thus, if $B$ is isotropic then $B = B'\oplus H$ for some form $B'$ (see Proposition 2.25 in \cite{gerstein}) and, consequently, $B=B' \in W(\mathbb F)$. If $B'$ itself is isotropic, there is a further decomposition $B'= B''\oplus H$ and again $B'=B'' \in W(\mathbb F)$. This process ends after a finite number of steps giving us a decomposition of the original form $B$, called the {\em Witt decomposition},  as 
$$B = B_0 \oplus H_1 \oplus ... \oplus H_n$$
where $B_0$ is anisotropic (but possibly zero) and each of $H_1,...,H_n$ (with $n$ also possibly zero) is a hyperbolic form. While this decomposition is not unique, the integer $n$ and the isomorphism type of $B_0$ are uniquely determined by $B$ (see Section 2.5 in \cite{gerstein}). With this in mind, we can define $W(\mathbb F)$ (as an Abelian group) alternatively as the set of equivalence classes of $\mathcal B_\mathbb F/\hspace{-2mm}\sim$ with the operation $\oplus$, where $\sim$ is defined as 
$$ B_1\sim B_2 \quad \Leftrightarrow \quad \text{The anisotropic parts of $B_1$ and $B_2$ are isomorphic.}$$
In this description it is easy to see that the inverse of $(V,B)$ in $W(\mathbb F)$ is the form $(V,-B)$ since $(V\oplus V, B+(-B))$ is a direct sum of hyperbolic forms. 
\vskip3mm

We now turn to examining some concrete Witt rings, including the case of $\mathbb F = \qq$. For a prime integer $p$, let $\mathbb Z_p$ denote the finite field $\mathbb Z_p = \{0,1,...,p-1\}$ of characteristic $p$. 
The Witt rings $W(\mathbb Z_p)$ are well understood as should be evident from the next theorem (which can be found in Section 2.8 in \cite{gerstein}). 
\begin{theorem} \label{theorem-witt-rings-for-finite-fields}
Let $p$ be a prime integer. Then there are isomorphisms of Abelian groups 
$$W(\mathbb Z_p) \cong \left\{ 
\begin{array}{cl}
\mathbb Z_2 & \quad ; \quad p = 2 \cr
\mathbb Z_2\oplus \mathbb Z_2 & \quad ; \quad p \equiv 1 \, \,( \text{mod } 4) \cr
\mathbb Z_4 & \quad ; \quad p \equiv 3 \, \,( \text{mod } 4)
\end{array}
\right. 
$$
The generators of $\mathbb Z_2\cong W(\mathbb Z_2)$ and of $\mathbb Z_4 \cong W(\mathbb Z_p)$ with $p  \equiv 3 \, \,( \text{mod } 4)$, are given by $\langle 1 \rangle$, while the two copies of $\mathbb Z_2$ in $W(\mathbb Z_p)$ when $p   \equiv 1 \, \,( \text{mod } 4)$, are generated by $\langle 1 \rangle$ and $\langle a \rangle$ for any choice of $a\in \dotz _p- (\dotz _p)^2$. 
\end{theorem}
The reason for stating a separate theorem about Witt rings of finite fields is that they are instrumental in understanding the Witt ring of the rationals. The relation between the former to the latter is elucidated in the next key theorem (which can be found on page 88 of \cite{husemoller-milnor}). 
\begin{theorem} \label{theorem-about-wq}
There is an isomorphism of Abelian groups 
$$ \sigma\oplus \partial : W(\mathbb{Q}) \to \mathbb{Z} \oplus  \left(  \oplus _{p}  W(\mathbb Z_p)  \right)$$
where $\oplus _p$ is a sum over all prime integers $p$. The homomorphism $\sigma :W(\mathbb{Q})\to \mathbb{Z}$ is the signature function while $\partial : W(\mathbb{Q})\to  \oplus _{p}  W(\mathbb Z_p) $ is the direct sum of homomorphisms $\partial _p :W(\mathbb{Q}) \to W(\mathbb Z_p)$ 
described on generators of $W(\qq)$ as follows: Given a rational number $\lambda \ne 0$, write it as $\lambda = p ^\ell \cdot \beta $ where  $\ell$ is an integer and $\beta$ a rational number whose numerator and denominator are relatively prime to $p$. Then 
\begin{equation} \label{delp}
\partial _p (\langle p ^\ell \cdot \beta \rangle ) = \left\{ 
\begin{array}{cl}
0 & \quad ; \quad \ell \mbox{ is even } \cr
\langle \beta \rangle & \quad ; \quad \ell \mbox{ is odd } 
\end{array}
\right. 
\end{equation}
\end{theorem}
The preceding theorem makes  is possible to determine precisely, and completely explicitly,  when two forms $B_1$ and $B_2$ over $\qq$ are equal in $W(\qq)$. Namely, if 
$$ B_1=\langle a_1 \rangle \oplus ... \oplus \langle a_n \rangle \quad \text{ and } \quad B_2=\langle b_1 \rangle \oplus ... \oplus \langle b_m \rangle $$
then $B_1=B_2 \in W(\qq)$ if and only if 
$$ \langle a_1 \rangle \oplus ... \oplus \langle a_n \rangle \oplus \langle - b_1 \rangle \oplus ... \oplus \langle - b_m \rangle = 0 \in W(\qq) $$
This latter equation in turn holds if and only if $\sigma$ and each $\partial _p$ map its left-hand side to zero. Here is an example illustrating Theorem \ref{theorem-about-wq}.  
\begin{example}
Let $B$ be the bilinear form on $\qq^4$ given by $B=\langle -\frac{23}{9}\rangle \oplus \langle 7 \rangle \oplus \langle -\frac{3}{5}\rangle \oplus \langle 49\rangle $. Note that the only primes $p$ for which $\partial _p B$ can be nonzero, are $p = 3, 5,7, 23$.  For these choices of $p$, we obtain 
\begin{align} \nonumber
\partial _3 B & = \langle \textstyle -\frac{1}{5} \rangle = \langle -5 \rangle = \langle 1 \rangle  \in W(\mathbb Z_3) \cr
\partial _5 B & = \langle -3 \rangle = \langle 2 \rangle  \in W(\mathbb Z_5)  \quad \text{ and } \quad 2\in \dotz_5 - (\dotz_5)^2 \cr 
\partial _7 B & = \langle 1 \rangle  \in W(\mathbb Z_7) \cr
\partial _{23} B & = \langle\textstyle -\frac{1}{9} \rangle = \langle -1 \rangle = \langle1 \rangle\oplus \langle1 \rangle\oplus\langle1 \rangle   \in W(\mathbb Z_{23})
\end{align}
Since $\sigma (B) = 0$, it follows that $B$ is torsion of order $4$ in $W(\qq)$. 
\end{example}

\section{The rational Witt class of a knot under a crossing change} \label{section-phi-of-k-under-a-crossing-change}
In this section we make precise the definition of $\varphi (K)$ -- the rational Witt class of a knot $K\subset S^3$ (see Definition \ref{definition-of-rational-witt-class-of-a-knot}). We then examine how $\varphi (K)$ is altered when $K$ undergoes a crossing change (Theorem \ref{how-things-are-after-a-crossing-change}). 
\vskip1mm

Let $K$ be an oriented knot in $S^3$ and let $\Sigma$ be a Seifert surface of $K$ whose orientation is compatible with that of $K$. We shall view the orientation of $\Sigma$ as being given by a normal and nowhere vanishing vector field $\vec{n}$ on $\Sigma$. The {\em linking form} or {\em Seifert form} on $H_1(\Sigma ; \mathbb Z)$ is the bilinear form $\ell k :H_1(\Sigma , \mathbb Z) \times H_1(\Sigma , \mathbb Z) \to \mathbb Z$ given by 
$$ \ell k (\alpha , \beta) = \text{ linking number of $\alpha$ with $\beta ^+$} $$
Here we view $\alpha$ and $\beta$ both as curves on $\Sigma$ and $\beta ^+$ is a small push-off of $\beta$ from $\Sigma$ in the direction of $\vec{n}$. Thus $\alpha$ and $\beta ^+$ are disjoint curves in $S^3$ and their linking number is $\frac{1}{2} \sum _{p} \eps(p)$ where $p$ ranges over the double points of any regular projection of $\alpha \sqcup \beta ^+$ and where $\eps(p)=1$ if $p$ is a positive crossing and $\eps(p) = -1$ is $p$ is a negative crossing (see Figure \ref{pic1} for the definition of positive/negative crossings). We extend $\ell k$ linearly to a form, of the same name, from $H_1(\Sigma ; \mathbb Q) \times H_1(\Sigma ; \mathbb Q)$ to $\qq$, and let $B_K : H_1(\Sigma ; \mathbb Q) \times  H_1(\Sigma ; \mathbb Q)  \to \mathbb Q$ be the bilinear, symmetric and non-degenerate form $B_K(\alpha, \beta ) =  \ell k (\alpha, \beta) + \ell k (\beta, \alpha)$. 
\begin{definition} \label{definition-of-rational-witt-class-of-a-knot}
With the notation as in the preceding paragraph, the rational Witt class $\varphi (K)$ of a knot $K\subset S^3$ is the element of the rational Witt ring $W(\qq)$ given by $(H_1(\Sigma ;\qq),B_K)$ for any choice of an oriented Seifert surface $\Sigma$ of $K$. 
\end{definition}
The fact that $\varphi (K)$ is well defined, i.e. independent of the choice of $\Sigma$, follows from work of Levine \cite{levine1, levine2} but can also be easily verified directly. Namely, any two oriented Seifert surfaces $\Sigma$ and $\Sigma '$ of the same knot $K$, differ from one another by a sequence of $1$-handle attachments/detachments. These operations change the associated bilinear forms by adding/subtracting a hyperbolic summand and thus do not affect their rational Witt classes. 

\begin{figure}[htb!] 
\centering
\includegraphics[width=14cm]{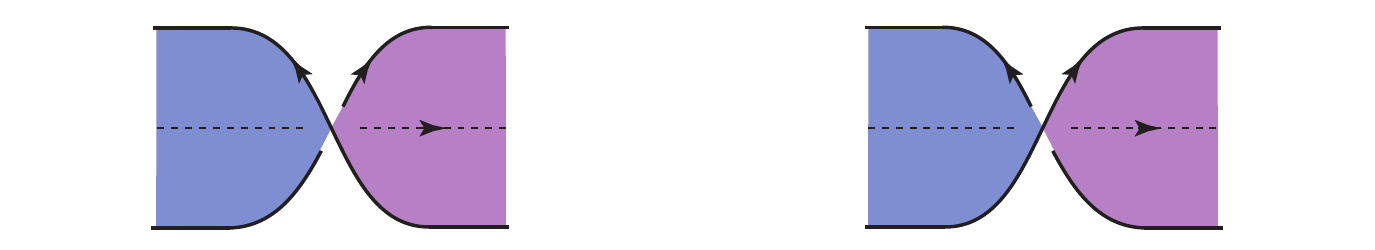}
\put(-305,48){$c$}
\put(-99,48){$c$}
\put(-270,20){$e_-^{2g}$}
\put(-64,20){$e_+^{2g}$}
\put(-350,45){$\Sigma_-$}
\put(-144,45){$\Sigma _+$}
\put(-350,-5){$K_-$}
\put(-144,-5){$K_+$}
\put(-310,-20){$(a)$}
\put(-105,-20){$(b)$}
\caption{Changing the negative crossing $c$ in $K_-$ to a positive one in $K_+$. The shaded areas indicate the Seifert surfaces $\Sigma _-$ and $\Sigma _+$.}  \label{pic2}
\end{figure}
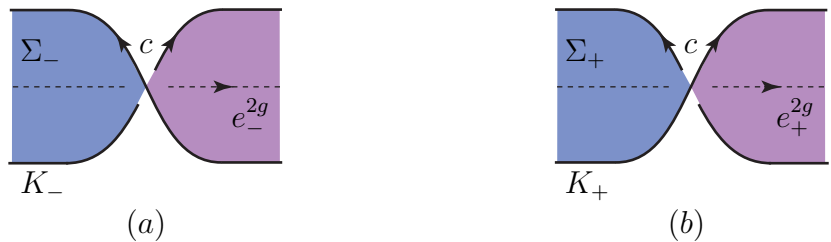

\vskip3mm
We now turn to exploring how $\varphi (K)$ changes when $K$ is altered by a single crossing switch. For concreteness sake, we take the crossing change to be a positive one, cf. Figure \ref{pic1}. Let $K_-$ be an oriented knot, let $c$ be a negative crossing in some projection of $K_-$ and let $K_+$ be the knot obtained from $K_-$ by switching the distinguished crossing $c$ as in Figure \ref{pic2}.

\begin{figure}[htb!] 
\centering
\includegraphics[width=14cm]{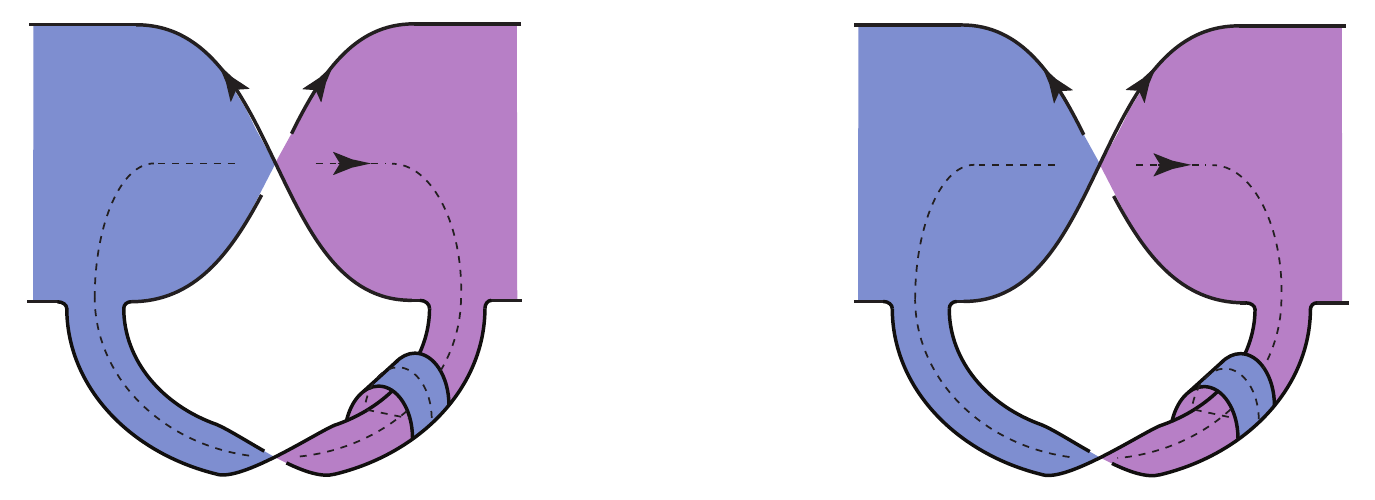}
\put(-370,100){$e^{2g}_-$}
\put(-133,100){$e^{2g}_+$}
\caption{One can always adjust the initial choices of Seifert surfaces $\Sigma _\pm$ by stabilizing them in a neighborhood of the distinguished crossing $c$, enabling one to find a preferred basis $\{e^1_\pm ,...,e^{2g}_\pm\}$ for $H_1(\Sigma _\pm ; \mathbb Z)$ with $e_-^i=e_+^i$ for $i=1,...,2g-1$. If initially a curve $e^i_\pm$ with $i\le 2g-1$ passes through the crossing $c$, we simply replace it by $e^i_\pm - e^{2g}_\pm$ (or by $e^i_\pm + e^{2g}_\pm$ depending on the orientations of the curves). }  \label{pic3}
\end{figure}

Pick oriented Seifert surfaces $\Sigma_\pm $ for $K_\pm $ so that $\Sigma _-$ and $\Sigma _+$ are identical safe in a neighborhood of the crossing $c$ where they differ as in Figure \ref{pic2}. For the purpose of comparing $\varphi (K_-)$ to $\varphi (K_+)$, it will prove advantageous to pick bases $\{e_\pm ^1,...,e_\pm^{2g}\}$ of $H_1(\Sigma_\pm;\mathbb{Z})$ with $e_-^i = e_+^i$ for $i=1,...,2g-1$, with $e_\pm^{2g}$ near $c$ as indicated in Figure \ref{pic3} and, additionally,  such that none of $e_\pm ^1, ..., e_\pm ^{2g-1}$ pass through the crossing $c$. Such bases can always be chosen though one may have to revise the initial choice of the Seifert surfaces $\Sigma _\pm$. Figure \ref{pic3} shows how to do this by a simple stabilization argument supported in a neighborhood of the distinguished crossing $c$. 

Let $\ell k_\pm : H_1(\Sigma _\pm;\mathbb Q) \times H_1(\Sigma _\pm;\mathbb Q)\to \mathbb Q$ be the linking pairings associated to $\Sigma _\pm$. Note that our choice of bases implies
\begin{align} \nonumber
\ell k_- (e^i_- , e^j_-) & = \ell k_+ (e^i_+ , e^j_+) \quad \quad \quad \quad \forall \, \,  (i,j) \ne (2g, 2g)  \cr
\ell k_- (e^{2g}_- , e^{2g}_-) & = \ell k_+ (e^{2g}_+ , e^{2g}_+) +1
\end{align}
Let $V_\pm$ be the $(2g)\times (2g)$ integral matrices representing the linking pairings $\ell k_\pm  $ with respect to the bases $\{e^1_\pm , ... , e^{2g}_\pm \}$, so that $\varphi (K_\pm ) = (\qq^{2g},V_\pm +V^\tau _\pm) \in W(\qq)$. We note that while the determinant of $\varphi (K)$ is only well defined as an element of $\dot\qq/\dot\qq^2$ (rather than as a rational number), the determinant of $V_\pm + V_\pm ^\tau$ agrees with the determinant of the knots $K_\pm$. 

To be able to compare $\varphi (K_-)$ to $\varphi (K_+)$, we shall express each as a sum of $1$-dimensional forms by diagonalizing $V_\pm +V_\pm ^\tau$. We accomplish this by changing our preferred bases $\{e_\pm ^1,...,e_\pm ^{2g}\}$ to new bases $\{f_\pm^1, ..., f_\pm ^{2g}\}$ via, essentially, the Gram-Schmidt algorithm. For simplicity of notation, we shall write $\langle v, w\rangle_\pm $ or simply $\langle v, w\rangle$ for $(\ell k _\pm + \ell k_\pm ^\tau)(v,w)$. 

With this in mind, we define the vectors $f^i_\pm$ as 
$$f^1_\pm = e^1_\pm \quad \quad \text{ and } \quad \quad f_\pm ^i = e_\pm ^i - \sum _{j=1}^{i-1} \frac{\langle e_\pm ^i, f_\pm ^j\rangle }{\langle f_\pm ^j, f_\pm ^j\rangle } f_\pm ^j \quad \text{ for } i\ge 2$$
These definitions may be ill posed since some of the numbers $\langle f^i_\pm , f^i_\pm\rangle$ could equal zero. To account for this, we divide our discussion into three separate cases. 

{\bf Case of $\langle f^i_\pm , f^i_\pm \rangle \ne 0$, $i=1,...,2g-1$. } For the moment, we assume that none of $\langle f^i_\pm , f^i_\pm \rangle$ vanishes. In this case we find that $\{f^1_\pm,...,f^{2g}_\pm\}$ are orthogonal bases for $H_1(\Sigma _\pm, \mathbb Q)$  where we think of $\ell k_\pm + \ell k ^\tau_\pm$ as giving us an inner product $\langle \cdot, \cdot \rangle = \lr_\pm$. Since $e^i_- = e^i_+$ for all $i=1,...,2g-1$, the same is true for the new basis elements: $f^i_- = f^i_+$ for $i=1,...,2g-1$. From this, and since $\langle e^i_-, e^j_-\rangle  = \langle e^i_+, e^j_+\rangle $ whenever $(i,j) \ne (2g,2g)$, we obtain 
$$ \langle f^i_-, f^j_- \rangle = \langle f^i_+ ,f^j_+ \rangle \quad \quad \forall \, \, (i,j) \ne (2g,2g)$$
On the other hand, for $i=j=2g$, we obtain 
\begin{align} \nonumber
& \langle f_+^{2g} , f_+^{2g}\rangle  = \left\langle  e_+^{2g} - \sum _{j=1}^{2g-1} \frac{\langle e_+^i, f_+^j\rangle }{\langle f_+^j, f_+^j\rangle } f_+^j \, , \, e_+^{2g} - \sum _{j=1}^{2g-1} \frac{\langle e_+^i, f_+^j\rangle }{\langle f_+^j, f_+^j\rangle } f_+^j             \right \rangle  \cr
& \quad = \langle e_+^{2g},e_+^{2g}\rangle - 2 \left\langle e_+^{2g} \, , \, \sum _{j=1}^{2g-1} \frac{\langle e_+^i, f_+^j\rangle }{\langle f_+^j, f_+^j\rangle } f_+^j \right\rangle  + 
\left\langle \sum _{j=1}^{2g-1} \frac{\langle e_+^i, f_+^j\rangle }{\langle f_+^j, f_+^j\rangle } f_+^j \, , \, \sum _{j=1}^{2g-1} \frac{\langle e_+^i, f_+^j\rangle }{\langle f_+^j, f_+^j\rangle } f_+^j             \right \rangle  \cr
& \quad =  \langle e_-^{2g},e_-^{2g}\rangle -2- 2 \left\langle e_-^{2g} \, , \, \sum _{j=1}^{2g-1} \frac{\langle e_-^i, f_-^j\rangle }{\langle f_-^j, f_-^j\rangle } f_-^j \right\rangle  + 
\left\langle \sum _{j=1}^{2g-1} \frac{\langle e_-^i, f_-^j\rangle }{\langle f_-^j, f_-^j\rangle } f_-^j \, , \, \sum _{j=1}^{2g-1} \frac{\langle e_-^i, f_-^j\rangle }{\langle f_-^j, f_-^j\rangle } f_-^j             \right \rangle  \cr
 & \quad  = \langle f_-^{2g}, f_-^{2g} \rangle -2 
\end{align}
Setting $\langle f^i_-,f^i_-\rangle = a_i$, we see that $\varphi (K_-)$ and $\varphi (K_+)$ take the forms 
\begin{align} \nonumber 
\varphi (K_-) & = \langle a_1 \rangle \oplus ... \oplus \langle a_{2g-1} \rangle \oplus \langle a_{2g} \rangle \cr
\varphi(K_+) & = \langle a_1 \rangle \oplus ... \oplus \langle a_{2g-1} \rangle \oplus \langle a_{2g} -2\rangle
\end{align}
Since we have taken care to make our basis change an orthogonal one, and since the determinant of $\varphi (K_\pm)$ when expressed in the bases $\{e^1_\pm, ..., e^{2g}_\pm \}$ agreed with $\det K_\pm$, we see that 
$$ \det K_- = | a_1 \cdot ... a_{2g-1}\cdot a_{2g}| \quad \text{and} \quad \det K_+ = |a_1 \cdot ... \cdot a_{2g-1}\cdot (a_{2g}-2)|$$
Given this, we can summarize our findings for the current special case by stating that, given our notation above, there exists a rational number $a\in \dot\qq$ (with $a=a_{2g}$ above) such that 
\begin{align} \label{forms-for-varphi-kplusminus}
\varphi(K_+)= \varphi(K_-) \oplus  \left\langle -\frac{1}{a} \right\rangle \oplus \langle a-2\rangle \quad \text{and} \quad \det K_+ = \det K_- \cdot \left| \frac{a-2}{a} \right|
\end{align}
As we shall see, this statement remains true in the two subsequent cases as well. 

{\bf Case of $\langle f^m_\pm, f^m_\pm \rangle =0$ for some $m\le 2g-2$. } This case has been addressed in Theorem 4.3 from \cite{jabuka1}. It is shown there that by passing to another basis, one can split off a hyperbolic summand from $(\ell k_\pm  + \ell k^\tau _\pm, H_1(\Sigma _\pm ; \qq))$ thereby ridding oneself of an element of square zero. Here are the specifics: Suppose that $f^i_\pm$, $i=1,...,m$  have been defined as above and satisfy $\langle f^i_\pm , f^i_\pm \rangle \ne 0$ for $i=1,...,m-1$ and that additionally  $\langle f^m_\pm , f^m_\pm \rangle = 0$. We then define a new basis $\{f^1_\pm,...,f_\pm^{m-1},f_\pm^{m}, \mathfrak e^{m+1}_\pm,...,\mathfrak e^{2g}_\pm \}$ according to 
\begin{align} \label{reduction-of-hyperbolic-summand}
\mathfrak e_\pm ^{m+1} & = e_\pm ^{m+1} - \sum _{j=1}^{m-1} \frac{\langle e_\pm ^{m+1},f_\pm ^j \rangle }{\langle f_\pm ^{j},f_\pm ^{j}\rangle} e_\pm ^j   \cr
\mathfrak g_\pm ^{m+k} & = e_\pm ^{m+k} - \sum _{j=1}^{m-1} \frac{\langle e_\pm ^{m+k},f_\pm ^j\rangle }{\langle f_\pm ^j , f_\pm ^j \rangle } e_\pm ^j \cr
\mathfrak e_\pm ^{m+k} & = \mathfrak g _\pm ^{m+k} -  \frac{\langle \mathfrak g_\pm ^{m+k},f_\pm ^{m}\rangle }{\langle \mathfrak e_\pm ^{m+1},f_\pm ^m\rangle} \mathfrak e_\pm ^{m+1} - \cr
& \quad \quad \quad \quad 
-\frac{\langle \mathfrak g_\pm ^{m+k} , \mathfrak e_\pm ^{m+1}\rangle \cdot \langle \mathfrak e_\pm ^{m+1} , f_\pm ^m \rangle  - \langle \mathfrak g_\pm ^{m+k},f_\pm ^m\rangle\cdot \langle \mathfrak e_\pm ^{m+1},\mathfrak e_\pm ^{m+1}\rangle }{\langle \mathfrak e_\pm ^{m+1},f_\pm ^m\rangle \cdot \langle \mathfrak e_\pm ^{m+1},f_\pm ^m\rangle } f_\pm ^m
\end{align}
with the last two equations valid for $k\ge 2$. An explicit computation (addressed in the proof of Theorem 4.3 in \cite{jabuka1}) shows that this basis decomposes as 
$$\{f_\pm ^1, ..., f_\pm ^{m-1}\}\cup \{ f^m_\pm, \mathfrak e^m_\pm \} \cup \{ \mathfrak e^{m+1}, ..., \mathfrak e^{2g}_\pm \} $$
with the spans of each of the three sets perpendicular (with respect to $\lr_\pm$) to the spans of the other two. Additionally, the restriction of $\lr_\pm$ to the span of $\{  f^m_\pm, \mathfrak e^m_\pm \} $ is a hyperbolic form showing that 
$$(\ell k _\pm + \ell k ^\tau _\pm, H_1(\Sigma _\pm ; \qq)) = (\lr_\pm|_{Span}, Span) \in W(\qq)$$
where $Span$ is $Span(f_\pm ^1, ..., f_\pm ^{m-1}, \mathfrak e^{m+1}, ..., \mathfrak e^{2g}_\pm)$. However, in passing from $H_1(\Sigma _\pm , \qq)$ to $Span(f_\pm ^1, ..., f_\pm ^{m-1}, \mathfrak e^{m+1}, ..., \mathfrak e^{2g}_\pm)$ we have eliminated the square zero vector $f^m$. Note that $ \mathfrak e^{2g}_\pm = e^{2g}_\pm + \mathfrak a_\pm $ with $\mathfrak a_- = \mathfrak a_+$ so that 
$$\langle \mathfrak e ^{2g}_+,\mathfrak e ^{2g}_+ \rangle = \langle \mathfrak e ^{2g}_- ,\mathfrak e ^{2g}_- \rangle +2 $$
still holds, allowing us to re-derive \eqref{forms-for-varphi-kplusminus} just as before (after first eliminating additional vectors of square zero, if any). 

{\bf Case of $\langle f^{2g-1}_\pm, f^{2g-1}_\pm \rangle = 0$. } 
If  $\langle f^{2g-1}_\pm, f^{2g-1}_\pm \rangle = 0$ then  $\langle f^{2g-1}_\pm, e^{2g}_\pm \rangle \ne 0$ since otherwise the from $\ell k _\pm + \ell k ^\tau _\pm$ would be degenerate. Thus, we introduce the new 
basis element $f_\pm ^{2g}$ as 
$$ f_\pm ^{2g} = e_\pm ^{2g} - \sum _{i=1}^{2g-2} \frac{\langle e_\pm ^{2g} , f_\pm ^i \rangle}{\langle f_\pm ^i , f_\pm ^i \rangle} f_\pm ^i $$
achieving $\langle f_\pm ^{2g} , f^i_\pm \rangle =0$ for $i=1,...,2g-2$. This shows that the spans of the two subsets $\{f_\pm^1,...,f_\pm ^{2g-2}\} $ and $\{f^{2g-1}_\pm , f^{2g} _\pm \}$ of our new basis, are orthogonal (again, with respect to $\lr _\pm$). Moreover, since $\langle f^{2g-1}_\pm, f^{2g-1}_\pm \rangle = 0$, the vectors $\{f^{2g-1}_\pm, f^{2g}_\pm\}$ span a hyperbolic space showing that (with $Span_\pm = Span (f^1_\pm, ..., f^{2g-2}_\pm )$)
$$
\begin{array}{c} 
\varphi (K_-) = (\ell k _- + \ell k _-^\tau, H_1(\Sigma _- ;\qq)) = ( \lr_- |_{Span_-}, Span _-) = \quad \quad \quad  \quad \cr
\quad \quad  \quad \quad = ( \lr_+ |_{Span_+}, Span _+) = (\ell k _+ + \ell k _+^\tau, H_1(\Sigma _+ ;\qq)) =\varphi (K_+)
\end{array}
$$
Even though we obtain $\varphi (K_-) = \varphi (K_+)$ in this case, both $\varphi (K_-)$ and $\varphi (K_+)$ still adhere to the form from \eqref{forms-for-varphi-kplusminus}, using a little trick. Namely, we can add the form $ \langle 4 \rangle \oplus \langle \textstyle -\frac{1}{4} \rangle$ to $\varphi (K_-)$ (since the former is trivial in $W(\qq)$) and observe that 
$$ \langle 4 \rangle \oplus \langle \textstyle -\frac{1}{4} \rangle =  \langle 4 \rangle \oplus \langle \textstyle -\frac{1}{4} -2 \rangle $$
since $-\frac{1}{4} -2 = -3^2\cdot \frac{1}{4}$ and therefore $\langle \textstyle -\frac{1}{4} - 2\rangle = \langle \textstyle -\frac{1}{4} \rangle$. Thus, even in this current case, $\varphi (K_-)$ and $\varphi (K_+)$ satisfy formula  \eqref{forms-for-varphi-kplusminus}. 
%
%
\vskip3mm
\noindent We see that each of the three cases we just discussed, leads to formula \eqref{forms-for-varphi-kplusminus}, thus proving the next proposition. 
\begin{proposition} \label{how-things-are-after-a-crossing-change}
Let $K_-$ and $K_+$ be two oriented knots that have projections that are identical save near one crossing $c$ which is a negative one for $K_-$ and a positive one for $K_+$. Then there exists a rational numbers $a \in \dot\qq$ such that 
$$\varphi(K_+)= \varphi(K_-) \oplus  \left\langle -\frac{1}{a} \right\rangle \oplus \langle a-2\rangle \quad \text{and} \quad \det K_+ = \det K_- \cdot \left| \frac{a-2}{a} \right|$$
\end{proposition}
Note that the preceding theorem shows that the signatures of $K_-$ and $K_+$ are related as 
$$ \sigma (K_+) = 
\left\{
\begin{array}{ll}
\sigma (K_-) - 2  & \quad ; \quad 0<a<2 \cr
& \cr
\sigma (K_-) & \quad ; \quad a<0 \, \,  \text{  or  } \, \,  a>2 
\end{array}
\right.
$$
Said differently, a positive crossing change in a knot,  decreases its signature by either $0$ or $2$. This observation proves the bound \eqref{lower-bounds-on-uK} for the case of $\omega =-1$ (see Definition \ref{definition-of-tristram-levine-signature} below). 

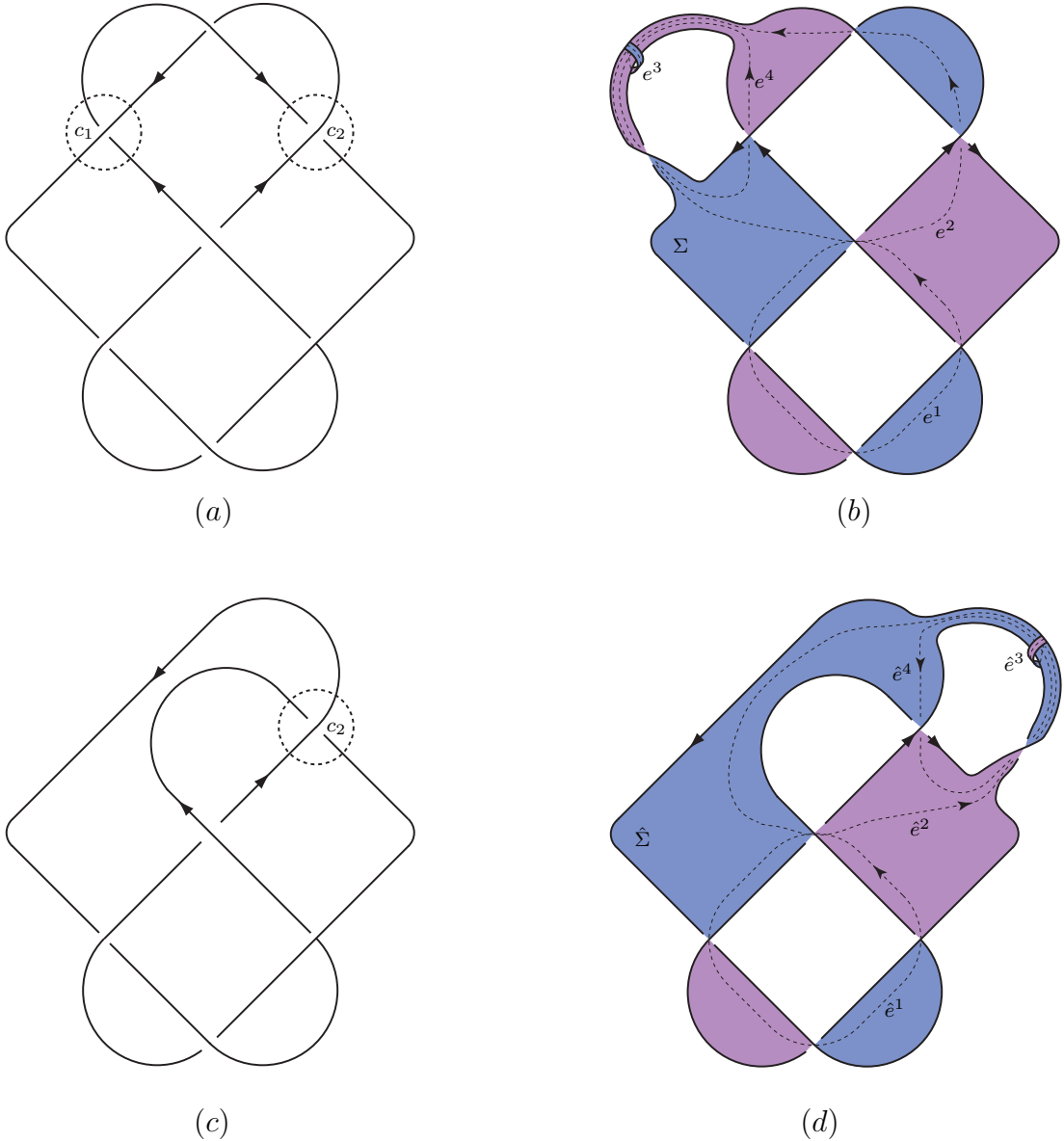
\begin{figure}[htb!] 
\centering
\includegraphics[width=15cm]{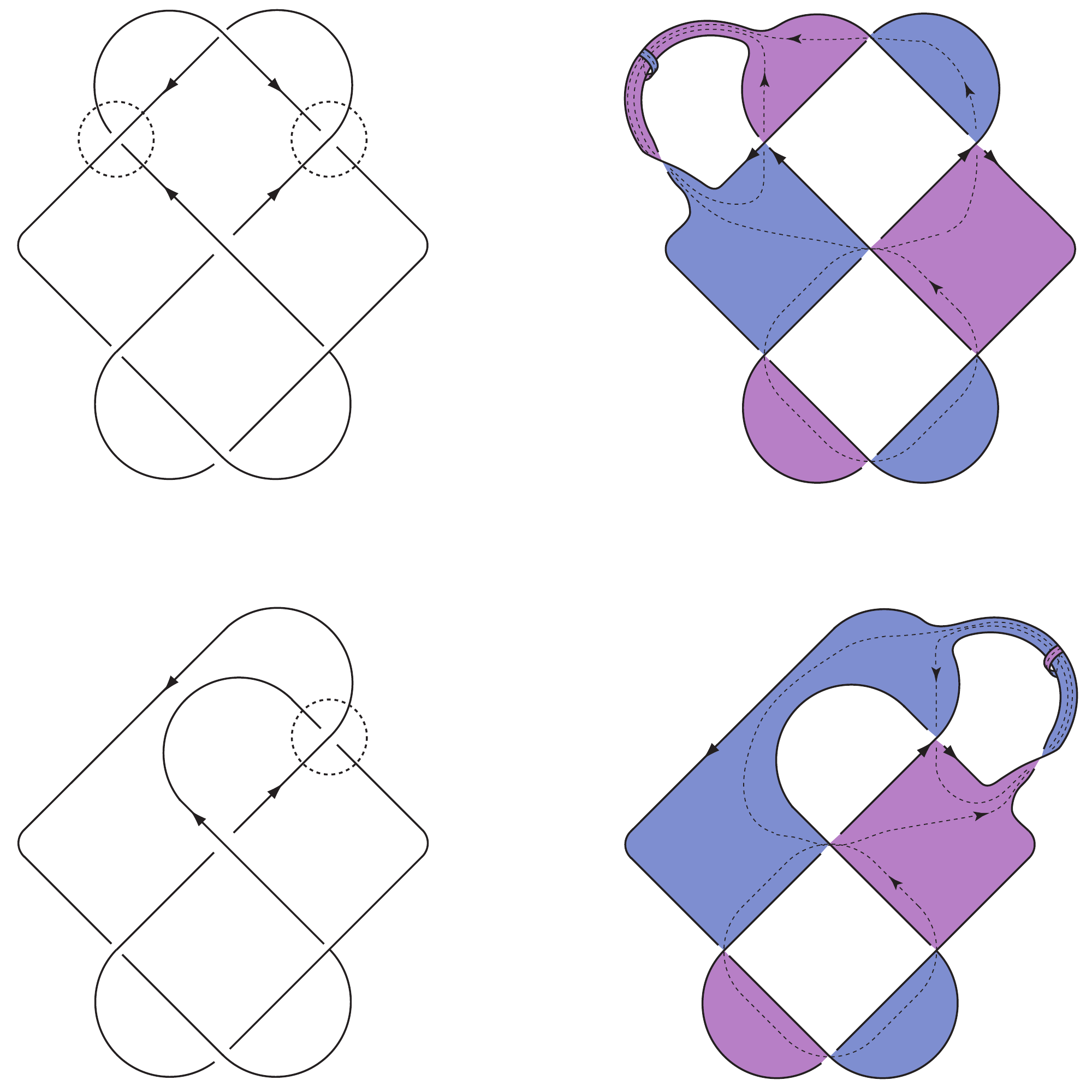}
\put(-393,371){\tiny $c_1$}
\put(-294,371){\tiny $c_2$}
\put(-294,137){\tiny $c_2$}
\put(-346,220){$(a)$}
\put(-94,220){$(b)$}
\put(-158,325){\tiny $\Sigma$}
\put(-173,92){\tiny $\hat \Sigma$}
\put(-60,258){\tiny $e^1$}
\put(-55,331){\tiny $e^2$}
\put(-170,393){\tiny $e^3$}
\put(-126,390){\tiny $e^4$}
\put(-346,-20){$(c)$}
\put(-108,-20){$(d)$}
\put(-75,25){\tiny $\hat e^1$}
\put(-65,96){\tiny $\hat e^2$}
\put(-28,161){\tiny $\hat e^3$}
\put(-72,157){\tiny $\hat e^4$}
\caption{{\bf (a)} The knot $7_4$ with two distinguished crossings $c_1$ and $c_2$. {\bf (b)} The Seifert surface $\Sigma$ and the basis $\{e^1, e^2, e^3, e^4\}$ of $H_1(\Sigma ; \mathbb {Z})$ used for analyzing how $\varphi (7_4)$ changes if the crossing $c_1$ is switched. The curve $e^3$ is oriented so that $\langle e^3,e^4\rangle =1$. {\bf (c)} The knot $7_4$ after the crossing $c_1$ has been switched, has undergone a simply isotopy showing the new knot to be $5_2$. Its remaining distinguished crossing $c_2$ is still indicated. {\bf (d)}  The Seifert surface $\hat \Sigma$ and the basis $\{ \hat e^1, \hat e^2, \hat e^3, \hat e^4\}$ of $H_1(\hat \Sigma ; \mathbb Z)$ used for describing how $\varphi (5_2)$ is affected by the change of the crossing $c_2$. The curve $\hat e^3$ is oriented so as to yield $\langle \hat e^3, \hat e^4\rangle =1 $.}  \label{pic5}
\end{figure}

These signature considerations along with the determinant formula from Proposition \ref{how-things-are-after-a-crossing-change}, allow for an explicit determination of $a$ from that same proposition. Namely, if $\sigma (K_+) = \sigma (K_-)$ then $\frac{a-2}{a}>0$ while if $\sigma (K_+) = \sigma (K_-)-2$ then $\frac{a-2}{a}<0$. From this one easily arrives at:
\begin{equation} \label{equation-solving-for-a}
\begin{array}{l}
\text{If } \sigma (K_+) = \sigma (K_-)  \text{ then } a=\frac{2 \det K_-}{\det K_- - \det K_+} \text{ and } a-2=\frac{2 \det K_+}{\det K_- - \det K_+}.  \cr \cr
\text{If } \sigma (K_+) = \sigma (K_-) - 2 \text{ then } a= \frac{2\det K_-}{\det K_- + \det K_+} \text{ and }  a-2 = -\frac{2\det K_+}{\det K_- + \det K_+}.
\end{array}
\end{equation}
Combining Proposition \ref{how-things-are-after-a-crossing-change} with the equations from \eqref{equation-solving-for-a}, provides a proof of the next theorem. 
\begin{theorem} \label{theorem-main1-before-r3}
Let $K_-$ and $K_+$ be two oriented knots that have projections that are identical safe near one crossing $c$ which is a negative one for $K_-$ and a positive one for $K_+$. Then their rational Witt classes are related as follows:  
$$\varphi(K_+)=
\left\{
\begin{array}{ll}
\varphi(K_-) \oplus  \left\langle -\frac{\det K_- - \det K_+ }{2 \det K-} \right\rangle \oplus \left\langle \frac{2 \det K_+}{\det K_- - \det K_+}\right\rangle &  \quad ; \quad \sigma (K_+) = \sigma (K_-) \cr  & \cr
\varphi(K_-) \oplus  \left\langle -\frac{\det K_- + \det K_+ }{2 \det K-} \right\rangle \oplus \left\langle -\frac{2 \det K_+}{\det K_- + \det K_+}\right\rangle &  \quad ; \quad \sigma (K_+) = \sigma (K_-) - 2 
\end{array}\right.
$$
\end{theorem}
Theorem \ref{main1} follows directly from the preceding theorem after an easy application of relation $(R3)$ from Theorem \ref{presentation-of-the-rational-witt-ring}. To see this, let us for brevity of notation write $\det K_- = b$ and $\det K_+ = c$. One then computes as
\begin{align} \nonumber
\textstyle \left\langle -\frac{\det K_- - \det K_+ }{2 \det K-} \right\rangle \oplus \left\langle \frac{2 \det K_+}{\det K_- - \det K_+}\right\rangle & = \textstyle \left\langle -\frac{b - c}{2 b} \right\rangle \oplus \left\langle \frac{2 c}{b -c}\right\rangle \cr
& = \langle -2b(b-c) \rangle \oplus \langle 2c(b-c)\rangle\quad  (\text{ now use $(R3)$} ) \cr
& = \langle -2(b-c)^2\rangle \oplus \langle 8bc(b-c)^4  \rangle \cr
& = \langle -2\rangle \oplus \langle 2bc\rangle \cr
& = \langle -2\rangle \oplus \left\langle \frac{2c}{b} \right\rangle \cr
& = \langle -2\rangle \oplus \left\langle\textstyle \frac{2\det K_+}{\det K_-} \right\rangle 
\end{align}
and more time as
\begin{align} \nonumber
\textstyle \left\langle -\frac{\det K_- + \det K_+ }{2 \det K-} \right\rangle \oplus \left\langle -\frac{2 \det K_+}{\det K_- + \det K_+}\right\rangle & = \textstyle \left\langle -\frac{b + c}{2 b} \right\rangle \oplus \left\langle -\frac{2 c}{b +c}\right\rangle \cr
& = \langle -2b(b+c) \rangle \oplus \langle -2c(b+c)\rangle\quad  (\text{ now use $(R3)$} ) \cr
& = \langle -2(b+c)^2\rangle \oplus \langle -8bc(b+c)^4  \rangle \cr
& = \langle -2\rangle \oplus \langle -2bc\rangle \cr
& = \langle -2\rangle \oplus \left\langle -\frac{2c}{b} \right\rangle \cr
& = \langle -2\rangle \oplus \left\langle \textstyle -\frac{2\det K_+}{\det K_-} \right\rangle 
\end{align}
With these in place, Theorem \ref{main1} follows. 
\vskip3mm
To illustrate our discussion thus far, we turn to an example. We shall consider the knot $K=\bar 7_4$ (the mirror image of the knot $7_4$) from Figure \ref{pic5}a and change the negative crossings $c_1$ and $c_2$ indicated in that same figure.   

To facilitate the crossing change at $c_1$, we pick the Seifert surface $\Sigma _{-}$ and the basis $\{e^1_{-},e^2_-, e^3_-,e^4_{-}\}$ for $H_1(\Sigma _{-};\qq)$ as indicated in Figure \ref{pic5}b (where we have dropped the subscripts \lq\lq $-$\rq\rq from the notation).  Let $V_{-}$ be the matrix representing the linking form $\ell k _{-}$ with respect to this basis. An explicit computation of linking numbers then shows that 
$$ V_{-} + V_{-}^\tau = \left[
\begin{array}{rr|rr}
4 & -1 & 0 & 0   \cr
-1 & 2 & 1 & -1   \cr \hline
0 & 1 & 0 & 1   \cr
0 & -1 & 1 & 0 
\end{array}
\right]
$$
Diagonalizing this matrix using the Gram-Schmidt process yields $\varphi (\bar 7_4)$ as 
$$\varphi (\bar 7_4) = \langle 4\rangle \oplus \langle \textstyle \frac{7}{4} \rangle \oplus \langle \textstyle -\frac{4}{7} \rangle \oplus \langle \textstyle \frac{15}{4} \rangle$$
Since changing the crossing $c_1$ only affects the linking number of $e^4_-$ with $e^4_-$, the rational Witt class of the knot $L$ obtained after changing $c_1$ differs from that for $7_4$ by subtracting $2$ from the $\frac{15}{4}$--summand:
$$\varphi (L) =  \langle 4\rangle \oplus \langle \textstyle \frac{7}{4} \rangle \oplus \langle \textstyle -\frac{4}{7} \rangle \oplus \langle \textstyle \frac{15}{4} - 2  \rangle =  \langle 4\rangle \oplus \langle \textstyle \frac{7}{4} \rangle \oplus \langle \textstyle -\frac{4}{7} \rangle \oplus \langle \textstyle \frac{7}{4} \rangle$$
This shows that $\varphi (\bar 7_4) = \varphi (L) \oplus \langle -\frac{1}{a}\rangle \oplus \langle a-2\rangle$ with $a=\frac{15}{4}$, as claimed in Theorem \ref{theorem-main1-before-r3}. It is rather easy to verify that $L$ is the knot $\bar 5_2$.

Turning now to changing the crossing $c_2$, we pick a new Seifert surface $\hat \Sigma _-$ and basis $\{\hat e^1_-,...,\hat e^4_-\}$ as in Figure \ref{pic5}(d) (where again the subscripts are omitted). Let $\hat V_-$ be the matrix expressing the linking form $\widehat{\ell k} _-$ with respect to this basis. A quick computation of linking numbers yields
$$ \hat V_{-} + \hat V_{-}^\tau = \left[
\begin{array}{rr|rr}
4 & -1 & 0 & 0   \cr
-1 & 2 & 1 & 1   \cr \hline
0 & 1 & 0 & 1   \cr
0 & 1 & 1 & 2 
\end{array}
\right]
$$
Running the Gram-Schmidt procedure on the latter bilinear form yields the rational Witt class $\varphi (\bar 5_2)$: 
$$\varphi (\bar 5_2) = \langle 4\rangle \oplus \langle \textstyle \frac{7}{4}\rangle \oplus \langle -\textstyle \frac{4}{7}\rangle \oplus \langle \textstyle \frac{7}{4} \rangle$$
The knot $M$ obtained from $\bar 5_2$ by changing $c_2$ must have rational Witt class equal to the one gotten from $\varphi (\bar 5_2)$ by subtracting $2$ from the last $\frac{7}{4}$--summand:
$$\varphi (M) =  \langle 4\rangle \oplus \langle \textstyle \frac{7}{4}\rangle \oplus \langle -\textstyle \frac{4}{7}\rangle \oplus \langle \textstyle \frac{7}{4} -2 \rangle =  \langle 4\rangle \oplus \langle \textstyle \frac{7}{4}\rangle \oplus \langle -\textstyle \frac{4}{7}\rangle \oplus \langle \textstyle -\frac{1}{4} \rangle = 0$$
By inspection one finds that $M$ is in fact that unknot and so $\varphi (\text{unknot}) = 0\in W(\qq)$, as already mentioned in the introduction. We note that here too, in accordance with Theorem \ref{theorem-main1-before-r3}, we obtain the equation $\varphi (\bar 5_2) = \varphi (M)\oplus \langle -\frac{1}{a}\rangle \oplus \langle a-2\rangle$, this time with $a=\frac{7}{4}$. 
\vskip3mm
We finish this section with a brief discussion of Tristram-Levine signatures and indicate, with few details,  how inequality \eqref{lower-bounds-on-uK} follows easily from the discussion preceding Theorem \ref{how-things-are-after-a-crossing-change} (we note that the bound \eqref{lower-bounds-on-uK} is only used with $\omega =-1$ in our proofs). 
\begin{definition} \label{definition-of-tristram-levine-signature}
Let $K$ be an oriented knot and $\Sigma$ an oriented Seifert surface for $K$.  Given a complex number $\omega$ of unit modulus, consider the Hermitian form $A(\omega)$ on $H_1(\Sigma;\mathbb C)$ given by
$$A(\omega) =  \textstyle \frac{1-\omega}{2} \cdot  \ell k +\textstyle \frac{1-\omega^{-1}}{2}  \cdot \ell k ^\tau$$
The Tristram-Levine signature $\sigma _\omega (K)$  is defined as the signature $\sigma (A(\omega))$ of $A(\omega)$, provided the latter is non-sigular. If $A(\omega)$ is singular, then we set $\sigma _\omega (K) = \frac{1}{2} ( \sigma _{\omega _-}(K) + \sigma _{\omega _+}(K))$ where $\omega _\pm$ are points on $S^1$ on either side of $\omega$ and sufficiently close to it. 
\end{definition}
It is not hard to verify that $\sigma _\omega (K)$ is independent of the choice of $\Sigma$ and that $A(\omega)$ is singular if and only if $\omega$ is a root of the Alexander polynomial $\Delta_K(t)$ of $K$.  Note also that $\sigma _{-1}(K)$ agrees with the usual signature $\sigma (K)$.  

If $K_-$ and $K_+$ are two knots that only differ in a single crossing $c$ which is a negative one for $K_-$ and a positive one for $K_+$, then one can diagonalize the corresponding $A_\pm(\omega)$  (this time as Hermitian rather than symmetric forms)  much as was already done for the case of $\omega =-1$ in the discussion leading up to Theorem \ref{how-things-are-after-a-crossing-change}. The corresponding analogue of equation \eqref{forms-for-varphi-kplusminus} is that there exist rational numbers $a_1,...,a_{2g}\in \dot\qq$ such that the diagonalized $A_-(\omega)$ and $A_{+}(\omega)$ look as
$$\langle a_1\rangle \oplus ... \oplus \langle a_{2g} \rangle \quad \text{ and } \quad \langle a_1\rangle \oplus ... \oplus \langle a_{2g} +(1-\text{Re} ( \omega ))\rangle$$
respectively. Since for $\omega \in S^1$ one obtains $0\le 1-\text{Re}(\omega)\le 2$, we see that either $\sigma _\omega (K_+) = \sigma _\omega (K_-)$ or $\sigma _\omega (K_+) = \sigma _\omega (K_-) + 2$. Inequality  \eqref{lower-bounds-on-uK} is an easy consequence of this observation. 
\section{Proofs of Corollaries \ref{coro1}, \ref{coro2} and \ref{coro3} }  \label{section-the-proofs-of-the-theorems}
\subsection{Proof of Corollary \ref{coro1}} 
Suppose that $K$ is an oriented knot in $S^3$ that can be unknotted with a single crossing change. According to \eqref{lower-bounds-on-uK} (with $\omega=-1$) we must have $\sigma (K) = 0$ or $\sigma (K) = \pm 2$. Let $U$ denote the unknot and recall that $\varphi (U)=0\in W(\qq)$.

If $K$ can be unknotted with a negative crossing change, then according to Theorem \ref{main1} (with $K_-=U$ and $K_+=K$) we obtain
$$
\varphi (K) = \left\{
\begin{array}{ll}
\langle  2 \det K \rangle \oplus  \langle - 2 \rangle & \quad ; \quad  \sigma (K) = 0\cr   & \cr
\langle  -2 \det K \rangle \oplus  \langle - 2 \rangle & \quad ; \quad  \sigma (K) = -2
\end{array}
\right.
$$
On the other hand, if $K$ can be unknotted with a positive crossing change, then Theorem \ref{main1} (this time with $K_-=K$ and $K_+=U$) implies  
$$
0 = \left\{
\begin{array}{ll}
\varphi(K) \oplus \langle  \frac{2}{\det K} \rangle \oplus  \langle - 2 \rangle & \quad ; \quad  \sigma (K) = 0\cr   & \cr
\varphi(K) \oplus \langle  -\frac{2}{\det K} \rangle \oplus  \langle - 2 \rangle & \quad ; \quad  \sigma (K) = 2
\end{array}
\right.
$$
Solving each of these for $\varphi (K)$ yields
$$
\varphi(K) = \left\{
\begin{array}{ll}
\langle - 2\det K \rangle \oplus  \langle 2 \rangle & \quad ; \quad  \sigma (K) = 0\cr   & \cr
\langle  2\det K \rangle \oplus  \langle 2 \rangle & \quad ; \quad  \sigma (K) = 2
\end{array}
\right.
$$
This completes the proof of Corollary \ref{coro1}. 
\subsection{Proofs of Corollaries \ref{coro2} and \ref{coro3}} 
We start by considering Corollary \ref{coro2} first. Let $K$ be a knot with unknotting number $2$ and let $L$ be the knot obtained from $K$ after a single crossing change. We let $U$ again denote the unknot. We split our discussion according to the type of crossing changes involved in changing $K$ to $U$ via $L$.

\subsubsection{$K$ can be unknotted with two negative crossing changes.} If $K$ can be unknotted with two negative crossing changes, then Theorem \ref{main1} implies that   
\begin{equation} \label{equation-toward-coro1-aux1}
\begin{array}{rl} 
\varphi (K) & = \left\{
\begin{array}{cl}
\varphi (L) \oplus \left\langle \frac{2\det K}{\det L} \right\rangle \oplus  \langle -2\rangle & \quad ; \quad \sigma (K) = \sigma (L)   \cr  & \cr
\varphi (L) \oplus \left\langle -\frac{2\det K}{\det L} \right\rangle \oplus  \langle -2\rangle & \quad ; \quad \sigma (K) = \sigma (L) - 2 
\end{array} \right. \cr 
& \cr
\varphi (L) & = \left\{
\begin{array}{cl}
\langle 2\det L \rangle \oplus  \langle -2\rangle & \quad ; \quad \sigma (L) = 0   \cr  & \cr
\langle -2\det L \rangle \oplus  \langle -2\rangle & \quad ; \quad \sigma (L) =  - 2 
\end{array} \right.
\end{array} 
\end{equation}
From these, part (a) of Corollary \ref{coro2} follows. For example, if $\sigma (K) = -4$ then $\sigma (L) = -2$ so that \eqref{equation-toward-coro1-aux1} implies
$$\varphi (K) =\textstyle \left\langle -\frac{2\det K}{\det L} \right\rangle \oplus \left\langle -2\det L \right\rangle \oplus  \langle -2\rangle  \oplus  \langle -2\rangle $$
Note that $\left\langle -\frac{2\det K}{\det L} \right\rangle = \left\langle -2\det K\cdot\det L \right\rangle \in W(\qq)$ and that, according to relation $(R3)$ from Theorem \ref{presentation-of-the-rational-witt-ring}, the equality $\langle -2\rangle \oplus \langle -2 \rangle = \langle -1\rangle \oplus \langle -1 \rangle$ also holds in $W(\qq)$. 

If $\sigma (K) = -2$ then either $\sigma (L) = -2$ or $\sigma (L) = 0$. These two cases, again in conjunction with \eqref{equation-toward-coro1-aux1}, lead to 
$$ 
\varphi (K) = \left\{
\begin{array}{ll}
\left\langle \frac{2\det K}{\det L} \right\rangle \oplus \langle -2\det L \rangle \oplus  \langle -2\rangle \oplus   \langle -2\rangle & \quad ; \quad \sigma (L) = -2 \cr   & \cr
 \left\langle -\frac{2\det K}{\det L} \right\rangle \oplus  \langle 2\det L \rangle \oplus  \langle -2\rangle \oplus \langle -2\rangle & \quad ; \quad \sigma (L) = 0 
\end{array}
\right. 
$$

Finally, if $\sigma (K) = 0$ then $\sigma (L) = 0$ also so that \eqref{equation-toward-coro1-aux1} provides us with 
$$\varphi (K) = \textstyle \left\langle \frac{2\det K}{\det L} \right\rangle \oplus  \langle 2\det L \rangle \oplus  \langle -2\rangle \oplus \langle -2\rangle$$
%
\subsubsection{$K$ can be unknotted by two positive crossing changes.} 
If $K$ is knot with $\sigma (K) \le 0$ (the assumption used in Corollary \ref{coro2}) and can be unknotted by two positive crossing changes, then $\sigma (K)=0$ since a positive crossing change cannot increase the signature. Consequently, we also obtain $\sigma (L) = 0$. Theorem \ref{main1} now implies that 
$$\varphi (K) = \varphi(L) \oplus \langle -2 \det K \det L\rangle \oplus \langle 2 \rangle \quad  \text{ and } \quad \varphi (L) = \langle -2 \det L \rangle \oplus \langle 2 \rangle $$
showing that $\varphi (K) =  \langle -2 \det K \det L\rangle \oplus  \langle -2 \det L \rangle \oplus \langle 1 \rangle \oplus \langle 1 \rangle$ (since $\langle 2 \rangle \oplus \langle 2 \rangle = \langle 1 \rangle \oplus \langle 1 \rangle$). 
%
\subsubsection{$K$ can be unknotted with one positive and one negative crossing change.} 
In this subsection we distinguish further the cases of $\sigma (K) = -2$ and $\sigma (K) = 0$. 

Starting with $\sigma (K) = -2$, we note that by signature considerations, it follows that the negative crossing change has to increase the signature of $K$ to $0$ while the positive crossing change cannot alter it further. Regardless of whether $L$ is gotten from $K$ by the positive or the negative crossing change, Theorem \ref{main1} implies that 
$$\varphi (K) = \langle -2 \det K \det L\rangle \oplus \langle -2 \det L\rangle $$

If $\sigma (K) = 0$ then there are two possibilities, namely, either both the positive and the negative crossing change alter the signature,  or else, neither does. In both cases, Theorem \ref{main1} yields
$$\varphi (K) = \langle \pm 2 \det K \det L\rangle \oplus \langle \mp \det L\rangle$$
%
%
%
%
This exhausts all possibilities and completes the proof of Corollary \ref{coro2}. 
\vskip3mm
The proof of Corollary \ref{coro3} follows along the same lines as the proofs of Corollaries \ref{coro1} and \ref{coro2}. Thus, suppose that $K$ is a knot with $\sigma (K) = -2n$ for some $n\in \mathbb N$ and that $u(K) = n$. Let $L_i$ be the knot gotten from $K$ by changing $i-1$ of these $n$ crossings so that, for instance, $L_1=K$ while $L_{n+1}$ is the unknot.  Signature considerations dictate that all of these crossing changes by negative crossing changes (since positive crossing changes cannot increase the signature) and that therefore $\sigma (L_i) = -2(n-i+1)$. Using this observation, Theorem \ref{main1} implies that 
$$\varphi (L_{i}) = \varphi (L_{i+1}) \oplus \langle -\textstyle \frac{2 \det L_{i+1}}{\det L_{i}} \rangle \oplus \langle -2\rangle   \quad \quad \quad  i=1, ..., n$$
Adding these last $n$ equations immediately yields the result of Corollary \ref{coro3}. 
\section{Low crossing examples}  \label{section-about-the-low-crossing-knots}
The results of Corollaries \ref{theorem-results-for-11-crossing-knots} and \ref{theorem-results-for12-crossing-knots} are a direct consequence of applying Corollary \ref{coro1} to certain $11$ and $12$ crossing knots.  Our computations were aided by a {\sc Mathematica} computer code written by the author. The Seifert matrices for the various knots were taken from KnotInfo\cite{knotinfo2}, indeed, without the latter our calculations would have been substantially more time consuming. In the next two subsections, we list with full details, two sample computations.
\subsection{Knots with $11$ crossings}
In this section we consider the knot $K=11a_{16}$ as an example. This knot has signature zero and determinant  $105$ and so, in order to have unknotting number $1$, its rational Witt class (according to Corollary \ref{coro1}) must equal 
$$\varphi (11a_{16}) = \langle \pm 210 \rangle \oplus \langle \mp 2 \rangle $$
for at least one consistent choice of signs. 

To compute the actual rational Witt class of $11a_{16}$, we start with a Seifert form the latter. From KnotInfo \cite{knotinfo2}, one finds that the symmetrized linking form $\ell k + \ell k ^\tau$ of $11a_{16}$ is represented by the matrix 
$$ \ell k + \ell k ^\tau  = \left[ 
\begin{array}{rrrrrr}
-2 & -1 & -1 & 0 & 0 & 1 \cr
-1 & -2 & -1 & 0 & 0 & 1 \cr
-1 & -1 & 2 & 0 & 1 & -1 \cr
0 & 0 & 0 & -2 & 0 & 1 \cr
0 & 0 & 1 & 0 & 2 & -1 \cr
1 & 1 & -1 & 1 & -1 & 4 
\end{array}
\right]
$$
This matrix is then diagonalized using the Gram-Schmidt procedure (without having to split off hyperbolic summands) to give 
$$ \varphi (11a_{16}) = \langle -2 \rangle \oplus \langle -2 \rangle \oplus  \left\langle -\frac{3}{2}\right\rangle \oplus \left\langle \frac{8}{3} \right\rangle \oplus \left\langle \frac{13}{8} \right\rangle  \oplus \left\langle \frac{105}{26} \right\rangle $$
To compare the latter to $\langle \pm 210 \rangle \oplus \langle \mp2\rangle$, we apply the homomorphism $\partial _5 :W(\qq) \to W(\mathbb Z_5)$  to all three forms: 
\begin{align}\nonumber
\partial _5 (\varphi (11a_{16})) & = \langle 21\cdot 26\rangle =  \langle 1\rangle  \cr
\partial _5( \langle \pm 210\rangle \oplus \langle \mp 2 \rangle  ) & = \langle \pm 42 \rangle = \langle 3 \rangle 
\end{align}
Since $1\in (\dotz_5)^2$ but $3\in \dotz _5 - (\dotz_5)^2$, the forms $\langle 1 \rangle$ and $\langle 3 \rangle$ are distinct forms in $W(\mathbb Z _5)$ (see Theorem \ref{theorem-witt-rings-for-finite-fields}). Accordingly, $\varphi (11a_{16})$ cannot equal $\langle 210\rangle \oplus \langle -2\rangle$ nor $\langle -210\rangle \oplus \langle 2\rangle$ in $W(\qq)$ and therefore $u(11a_{16})\ge 2$. Since an explicit unknotting of $11a_{16}$ with two crossing changes is easily found, we arrive at $u(11a_{16}) = 2$ as claimed in Corollary \ref{theorem-results-for-11-crossing-knots}. 
\subsection{Knots with $12$ crossings}
As an example among $12$ crossings knots, we single out the non-alternating knot $K=12n_{33}$. This knot has signature $-2$ and determinant $123$. If we had $u(12n_{33})=1$, then according to Corollary \ref{coro1}, its rational Witt class would have to equal $\langle -246\rangle \oplus \langle -2\rangle$. The actual rational Witt class of $12n_{33}$ is again computed by starting with the matrix representing $\ell k + \ell k ^\tau$ which one finds on KnotInfo \cite{knotinfo2} to be
$$\ell k +\ell k ^\tau  = \left[
\begin{array}{rrrrrrrr}
2 & -1 & 0 & -1 & 0 & 0 & -1 & -1 \cr
-1 & -2 & 0 & -1 & 0 & 0 & -1 & -1 \cr
0 & 0 & -2 & -1 & 0 & 0 & -1 &- 1\cr
-1 & -1 & -1 & -2 & -1 & -1 & -1 & -1 \cr
0 & 0 & 0 & -1 & 2 & 1 & -1 & -1 \cr
0 & 0 & 0 & -1 & 1 & 2 & 0 & 0 \cr
-1 & -1 & -1 & -1 & -1& 0 & -2 & -1 \cr  
-1 & -1 & -1 & -1 & -1 & 0 & -1 & -2
\end{array}
\right]
$$
This, when diagonalized with the Gram-Schmidt algorithm, yields (after a cancellation of $\langle 2 \rangle \oplus \langle -2 \rangle$) 
$$\varphi (12n_{33}) = \left\langle -\frac{5}{2} \right\rangle \oplus \left\langle -\frac{11}{10} \right\rangle \oplus \left\langle \frac{2}{11} \right\rangle \oplus \left\langle \frac{53}{2}\right\rangle \oplus \left\langle -\frac{22}{53} \right\rangle \oplus \left\langle -\frac{123}{22}\right\rangle$$
Applying $\partial _{41}$ to these gives 
\begin{align} \nonumber
\partial _{41}(\langle -246\rangle \oplus \langle -2 \rangle ) & = \langle -6 \rangle   = \langle 35 \rangle \cr
\partial _{41}(\varphi (12n_{33}) )& = \langle -3 \cdot 22 \rangle  = \langle 16\rangle 
\end{align}
Since $16$ is a square in $\mathbb Z_{41}$ while $35$ isn't, the two forms $\varphi (12n_{33})$ and $\langle -246\rangle \oplus \langle -2 \rangle $ cannot be equal in $W(\qq)$ (cf. Theorem \ref{theorem-witt-rings-for-finite-fields}) and so $u(12n_{33}) \ge 2$ as stated in Corollary \ref{theorem-results-for12-crossing-knots}. 
\subsection{Proof of Corollary \ref{corollary-for-unknotting-number-two}}
In Corollary \ref{corollary-for-unknotting-number-two} we examined the knot $10_{47}$. For convenience, we work here with the knot $K=\overline{10}_{47}$ (the mirror image of the knote $10_{47}$) which has signature $-4$ and determinant $41$. If we had $u(K)=2$, Corollary \ref{coro2} implies that the rational Witt class of $K$ must be equal to 
\begin{equation} \label{the-first-form}
\varphi (\overline{10}_{47}) = \langle -82d\rangle \oplus \langle -2d \rangle \oplus \langle -1 \rangle \oplus \langle -1 \rangle 
\end{equation}
where $d=\det L$ is the determinant of the knot $L$ obtained from $K$ after only one crossing change. The signature of $L$ is clearly $-2$. On the other hand, the actual rational Witt class of $\overline{10}_{47}$ can be computed to be
\begin{equation} \label{the-second-form}
\varphi (\overline{10}_{47}) =  \langle 2\rangle \oplus \langle \textstyle \frac{3}{2}\rangle \oplus \langle \textstyle -\frac{8}{3} \rangle \oplus  \langle \textstyle -\frac{13}{8} \rangle \oplus  \langle \textstyle -\frac{18}{13}  \rangle \oplus  \langle \textstyle -\frac{27}{18}  \rangle \oplus  \langle \textstyle -\frac{34}{27}  \rangle \oplus  \langle \textstyle -\frac{41}{34}\rangle 
\end{equation}
If one assumes that $L$ is a knot with $9$ or fewer crossings, then $d$ takes on odd values from $1$ to $75$ (as can be seen by examining the tables at KnotInfo \cite{knotinfo2}). One then asks, for which values of $d$ in that range, are the rational Witt classes from \eqref{the-first-form} and \eqref{the-second-form} equal. Using Mathematica, one finds the answer is for 
$$d=3, 7, 11, 15, 19, 27, 35, 47, 55, 63, 67, 71, 75$$
Finally, the only knots $L$ with $9$ or fewer crossings, with unknotting number $1$, with signature $-2$ and with determinant given by one of the $d$'s from the previous line, are the knots 
$$L = 3_1, 5_2,6_2, 7_2, 7_6, 8_{11}, 8_{21}, 9_2, 9_{12}, \bar 9_{26}, \bar 9_{39},  \bar 9_{42}$$
as claimed in Corollary \ref{corollary-for-unknotting-number-two}.

\section{Pretzel knots} \label{section-on-pretzel-knots}
The rational Witt classes and the signatures of pretzel knots $P(p_1,...,p_n)$ have been completely determined in \cite{jabuka1}. The pertinent statements are contained in Theorems 1.2 -- 1.4 and in Theorem 1.18 in \cite{jabuka1}. 
For the reader's convenience, we provide below those results relevant to Section \ref{introductory-section-on-pretzel-knots} of the present article.  
\subsection{Case $3$-stranded pretzel knots }  Given a pretzel knot $P(p_1,p_2,p_3)$ with $p_1,p_2$ odd and $p_3\ne 0$ even, its rational Witt class and signature are given by (courtesy of \cite{jabuka1}): 
\begin{align} \nonumber
\varphi (P(p_1,p_2,p_3)) & = \bigoplus _{i=1}^2 \left( \oplus _{k=1}^{|p_i|-1}  \langle -\eps _k k(k+1) \rangle \right)  \oplus \langle -\textstyle \frac{p_1+p_2}{p_1p_2}\rangle \oplus \langle \frac{\det P(p_1,p_2,p_3)}{p_1+p_2}\rangle \cr
\sigma (P(p_1,p_2,p_3)) & = (\eps_1+\eps_2) - (p_1+p_2) -Sign (\textstyle \frac{p_1+p_2}{p_1p_2} ) + Sign(\textstyle\frac{det P(p_1,p_2,p_3)}{p_1+p_2})
\end{align}
In the above, $\eps _k = Sign (p_k)$ and $\det P(p_1,p_2,p_3) = p_1p_2+p_1p_3+p_2p_3$ is a signed version of the determinant used in \cite{jabuka1}. We have implicitely assumed that $p_1+p_2\ne 0$ for if $p_1+p_2=0$, then $\varphi ( P(p_1,-p_1,p_3))=0$ and $\sigma (P(p_1,-p_1,p_3))=0$. 

The knots considered in Corollary \ref{pretzel1} make the choices $p_1\ge 7$,  $p_2=4-p_1$ and $p_3>-\frac{p_1(4-p_1)}{4}$, the latter condition ensuring that $\det P(p_1,4-p_1,p_3)>0$. The rational Witt class and signature  of $P(p_1,4-p_1,p_3)$ then become 
\begin{align} \label{form-for-varphi-of-pp14-p1p3-first-version}
\varphi  (P(p_1,4-p_1,p_3)) & =  \bigoplus _{i=1}^2 \left( \oplus _{k=1}^{|p_i|-1}  \langle -\eps _k k(k+1) \rangle \right)  \oplus \langle -p_1(4-p_1)\rangle \oplus \langle 4p_3+p_1(4-p_1) \rangle\cr
\sigma(P(p_1,4-p_1,p_3)) & = -2
\end{align}
We shall simplify the expression for $\varphi (P(p_1,4-p_1,p_3))$ by using the next lemma. 
\begin{lemma} \label{auxxlemma1}
For any $n\ge2$ and for $\eps \in \{\pm 1\}$, the equality 
$$ \oplus _{k=1}^{n-1} \langle -\eps \cdot  k(k+1)\rangle = \langle \eps \cdot n\rangle \oplus \left(\bigoplus _{i=1}^{n} \langle - \eps \rangle \right)$$
holds in $W(\qq)$. 
\end{lemma}
\begin{proof}
The claim of the lemma follows easily from an induction argument. When $n=2$, the equality 
$$ \langle -\eps \cdot 2\rangle = \langle \eps \cdot 2 \rangle \oplus \langle -\eps \rangle \oplus \langle -\eps \rangle $$
follows from an application of relation $(R3)$ from Theorem \ref{presentation-of-the-rational-witt-ring}, by which $\langle -2\eps\rangle \oplus \langle -2\eps \rangle = \langle -\eps \rangle \oplus \langle -\eps \rangle$. Proceeding by induction, and using again $(R3)$, we compute:
\begin{align} \nonumber
 \oplus _{k=1}^{n-1} \langle -\eps  k(k+1)\rangle & =  \langle -\eps (n-1)n\rangle \oplus \left( \oplus _{k=1}^{n-2} \langle -\eps k(k+1)\rangle \right) \cr
& =  \langle -\eps (n-1)n\rangle \oplus \langle \eps (n-1) \rangle \oplus \left(\bigoplus _{i=1}^{n-1} \langle - \eps \rangle \right) \quad \text{(now use $(R3)$) } \cr
& = \langle -\eps\rangle \oplus \langle \eps n \rangle \oplus \left(\bigoplus _{i=1}^{n-1} \langle - \eps \rangle \right) \cr
& = \langle \eps n \rangle \oplus \left(\bigoplus _{i=1}^{n} \langle - \eps \rangle \right)
\end{align}
This proves the lemma.  
\end{proof}
Lemma \ref{auxxlemma1} allows us to substantially simplify the expression for $\varphi (P(p_1,4-p_1,p_3))$ from \eqref{form-for-varphi-of-pp14-p1p3-first-version}, to obtain
$$\varphi (P(p_1,4-p_1,p_3)) = \langle p_1\rangle \oplus \langle 4-p_1\rangle \oplus \langle -p_1(4-p_1)\rangle \oplus \langle 4p_3+p_1(4-p_1) \rangle \oplus  \left(\bigoplus _{i=1}^{4} \langle -1 \rangle \right)$$
Applying relation $(R3)$ to the first two terms on the right-hand side above (and carrying the third term), yields an additional simplification:
$$
 \langle p_1\rangle \oplus \langle 4-p_1\rangle  \oplus \langle -p_1(4-p_1)\rangle   = \langle 4 \rangle \oplus \langle p_1(4-p_1)\rangle  \oplus \langle -p_1(4-p_1)\rangle = \langle 1 \rangle 
$$
With this last equation, we finally arrive at 
$$\varphi (P(p_1,4-p_1,p_3)) = \langle -1 \rangle \oplus \langle -2 \rangle \oplus \langle -2 \rangle \oplus \langle 4p_3+p_1(4-p_1) \rangle \oplus  \left(\bigoplus _{i=1}^{3} \langle -1 \rangle \right)$$
Corollary \ref{pretzel1} follows immediately from this and from Corollary \ref{coro1}. 

To explain the conclusions from Example \ref{example-one-for-pretzels-with-p3-even}, consider the knot $P(7,-3,r)$ with $r\ge 6$. Then $u(P(7,-3,r))=1$ forces the equality 
$$ \langle -1\rangle \oplus \langle -2 \rangle \oplus \langle 4r-21)\rangle = \langle -2(4r-21)\rangle $$
in $W(\qq)$. If there is a prime $p$ with $4r-21 = p^{2m+1}\cdot \beta$ and with $\gcd (\beta, p) =1$, the $\partial _p$ applied to the above, leads to the equality 
$$  \langle \beta \rangle = \langle -2 \beta \rangle \quad \text{ in } \quad W(\mathbb Z_p)$$
This latter equality can only be valid if $-2$ is a square in $\mathbb Z_p$. If $r=2k\cdot 7^{\ell +1}$ with $k,\ell \in \mathbb N$, then $4r-21 = 7(8k\cdot 7^\ell -3)$ so we can choose $p=7$ since, indeed, $-2$ is not a square in $\mathbb Z_7$. Example \ref{example-two-for-pretzels-with-p3-even} is analyzed similarly. 
\subsection{Case of $4$-stranded pretzel knots} 
Consider here a $4$-stranded pretzel knot $P(p_1,p_2,p_3,p_4)$ with $p_1,p_2,p_3$ odd and with $p_4\ne 0$ even. In this case, the rational Witt class and signature are given by (as proved in \cite{jabuka1}, with the use of Lemma \ref{auxxlemma1} for a slight simplification):
\begin{align} \nonumber
\varphi (P(p_1,p_2,p_3,p_4)) & = \bigoplus _{i=1}^4 \left(\langle p_i \rangle \oplus \bigoplus _{k=1}^{|p_i|} \langle - \eps _k \rangle     \right) \oplus \left\langle -\textstyle \frac{\det P(p_1,p_2,p_3,p_4)}{p_1p_2p_3p_4 }  \right\rangle  \cr
\sigma (P(p_1,p_2,p_3,p_4)) & = \sum _{i=1}^4 (\eps _i - p_i) - Sign (\,p_1p_2p_3p_4\, \det P(p_1,p_2,p_3,p_4)) 
\end{align}
where again $\eps _k = Sign (p_k)$ and where this time $\det P(p_1,p_2,p_3,p_4)=p_1p_2p_3+p_1p_2p_4+p_1p_3p_4+p_2p_3p_4$. 

Turning to the knot $P(p,p,p,-3p-1)$ considered in Corollary \ref{pretzel2}, we note that its determinant is given by 
$-p^2(8p+3)$. We remind the reader that $p$ is an odd, positive integer. When inserted into the above form of  the rational Witt class, this determinant formula provides us with 
\begin{align}\nonumber
\varphi(P(p,p,p,-3p-1)) & = \langle p\rangle \oplus  \langle p\rangle \oplus  \langle p\rangle \oplus \langle -3p-1\rangle \oplus \left\langle \textstyle - \frac{8p+3}{p(3p+1)}  \right\rangle  \oplus \langle 1 \rangle \cr
\sigma (P(p,p,p,-3p-1)) & = 2
\end{align}
Corollary \ref{pretzel2} follows directly from these formulas (with the use of Corollary \ref{coro1}). 

The validity of Example \ref{example-for-4-stranded-pretzels} is now easily deduced from Corollary \ref{pretzel2}. Namely, applying the homomorphism $\partial _{19}$ to $\varphi (P(p,p,p,-3p-1))$ with $p=2+(2k+1)\cdot 19^{\ell+1}$, yields 
$$\partial _{11} \left( \varphi (P(p,p,p,-3p-1) \right)= \langle 5 \rangle \in W(\mathbb Z_{19})$$
On the other hand, if we had $u(P(p,p,p,-3p-1))=1$, then Corollary \ref{coro1} would force the equality 
$\varphi (P(p,p,p,-3p-1)) = \langle 2 \rangle \oplus \langle 2(8p+3)\rangle $. However, $\partial _{19}$ applied to this last form gives (with $p=2+(2k+1)\cdot 19^{\ell+1}$)
$$ \partial _{19} ( \langle 2 \rangle \oplus \langle 2(8p+3)\rangle )= \langle 2 \rangle \in W(\mathbb Z_{19}) $$
Since $5$ is a square in $\mathbb Z_{19}$ but $2$ is not, we see that the equality $\varphi (P(p,p,p,-3p-1)) = \langle 2 \rangle \oplus \langle 2(8p+3)\rangle $ cannot be satisfied in $W(\qq)$. The conclusion of Example \ref{example-for-4-stranded-pretzels} follows. 
\subsection{Upward stabilizations} 
The notion of upward stabilization, used in Corollary \ref{pretzel3}, was introduced in Definition 1.6 from \cite{jabuka1}. As noted in \cite{jabuka1}, if $L$ is obtained by an upward stabilization from the pretzel knot $K$, then $\varphi (L) = \varphi (L)$ while $\det L = \det K \cdot \lambda ^2$ for some integer $\lambda$. These facts make Corollary \ref{pretzel3} evident. 
\section{Comparison with work of Lickorish \cite{Lickorish}} \label{section-work-of-lickorish}
This final section compares our work to that of R. Lickorish from \cite{Lickorish}. To set up the framework for comparison, we explore a few preliminaries first. 

Given a rational homology $3$-sphere $Y$, let 
$$\lambda : H_1(Y;\mathbb {Z}) \times H_1(Y;\mathbb {Z}) \to \qq/\mathbb Z$$
be its associated {\em linking pairing} defined as follows: Given two curves $\alpha, \beta \in H_1(Y;\mathbb Z)$, there exists a nonzero integer $n$ such that $n\cdot \alpha $ bounds a $2$-chain, say $\sigma$. With such a $\sigma$ in place, we set  
$$\lambda (\alpha, \beta) = \frac{\sigma \cdot \beta}{n}  \in \qq/\mathbb Z$$ 
where the numerator of the right hand side above is the usual intersection pairing between homology groups of complementary dimensions.  For our intentions, $Y$ will be the $2$-fold branched cover $\Sigma_K$ of $S^3$ with branching set a knot $K$.  

Lickorish proves the following lemma (Lemma 2 in \cite{Lickorish}):
\begin{lemma}[Lickorish \cite{Lickorish}] \label{lickorishs-lemma}
Let $Y$ be the $3$-manifold obtained by $p/q$-framed Dehn surgery on a knot $L\subset S^3$, with $p\ne 0$. Then $H_1(Y;\mathbb Z)$ is cyclic of order $|p|$, generated by a meridian $\mu$ and moreover
$$\lambda (\mu, \mu) = \frac{q}{p} \in \qq/\mathbb Z$$
\end{lemma}
This lemma establishes a bridge towards studying unknotting number $1$ knots since, if $u(K)=1$, then $\Sigma_K$ is obtained as $n/2$-surgery on some knot $L$ with $n$ an odd integer. This fact was already known to Montesinos \cite{Montesinos2, Montesinos} but a complete proof is also given by Lickorish in \cite{Lickorish}. Thus, according to Lemma \ref{lickorishs-lemma},  if $u(K) =1$, then $H_1(\Sigma_K;\mathbb{Z})$ is cyclic of order $|n|$ (and hence $n= \pm \det K$) and possesses a generator $\mu$ with $\lambda (\mu, \mu) = 2/n$.

If $K$ happens to be a $2$-bridge knot with $\Sigma _K = L(p,q)$, then $\Sigma _K$ is also gotten by $p/q$-framed surgery on the unknot and therefore, again according to Lemma \ref{lickorishs-lemma}, there has to be a generator $\mu '$ of $H_1(\Sigma _K;\mathbb Z)$ with $\lambda (\mu ', \mu ') = q/p$. Since $H_1(\Sigma _K,;\mathbb Z)$ is cyclic, there is an integer $t$ such that $\mu ' = t\cdot \mu$ and hence, by applying $\lambda$ to this, we also find that 
$$ \frac{q}{p} = t^2 \cdot \frac{2}{n}  \quad  \text{ in }  \quad \qq/\mathbb Z$$
This equation is re-captured in an equivalent format (since $n,p =\pm \det K$) by the next statement.
\begin{theorem}[Lickorish \cite{Lickorish}]   \label{lickorishs-theorem}
Let $K$ be a $2$-bridge knot with $2$-fold brached cover the lens space $L(p,q)$. If $u(K) = 1$, then the congruence 
$$ q \equiv \pm 2 t^2 \, \, \, (\text{mod } \det K) $$
must hold for some $t\in \mathbb Z$. 
\end{theorem}
Lickorish in \cite{Lickorish} goes on to apply this theorem to the knot $K=7_4$ for which $\Sigma _K = L(15,4)$. Thus, if $7_4$ had unknotting number $1$, there would have to be a solution $t\in \mathbb Z$ of the congruence 
$$ 4 \equiv \pm 2t^2\, \, \, (\text{mod } 15 ) $$
It is easy to see that there is no such $t$ showing that $u(7_4)>1$ and hence $u(7_4)=2$ (since an unknotting of $7_4$ with two crossing changes is easily found). 
\vskip3mm
Our methods equally well apply to the example $K=7_4$. Namely, its rational Witt class was computed in Section \ref{section-phi-of-k-under-a-crossing-change} as
$$\varphi (7_4) = \langle -4\rangle \oplus \langle -\textstyle\frac{15}{4} \rangle$$
If we had $u(7_4)=1$, then Corollary \ref{coro1} would force the equality $\langle -4\rangle \oplus \langle -\textstyle\frac{15}{4} \rangle = \langle -2\rangle \oplus \langle -\textstyle 2\cdot 15 \rangle$ in $W(\qq)$. This equality is easily seen to fail by applying $\partial _5$ to it, reducing it to $\langle 2\rangle  = \langle 1\rangle$ in $W(\mathbb Z_5)$ (which fails since $2$ is not a square in $\mathbb Z_5$). 
\vskip3mm
Lickorish's obstruction and our obstruction to the unknotting number being $1$, seem similar, at least in their origins, as both derive from an intersection pairing represented by the linking form of the knot. Yet, the two obstructions are not equivalent as the next example shows. 
\begin{example}
Consider the knot $K=8_8$ whose twofold branched cover is the lens space $L(25,9)$. If the unknotting number of $8_8$ were $1$, Theorem \ref{lickorishs-theorem} would guarantee a solution of the congruence $9\equiv \pm 2t^2 \, \, (\text{mod } 25)$. But there is no such solution since neither $2$ nor $23$ are squares mod $25$. It follows that $u(8_8) \ge 2$ (in fact, $u(8_8) = 2$). 

On the other hand, the obstruction to $u(8_8)=1$ from Corollary \ref{coro1}, just yields $\varphi (8_8) = 0$, an equality that is in fact valid. 
\end{example} 
Finding possible examples of two-bridge knots for which our obstruction provides positive results where Lickorish's doesn't, will have to wait until the rational Witt classes of two-bridge knots have been computed (a project currently in progress). There are, however, no such examples among two-bridge knots with $10$ or fewer crossings.


\end{document}